\documentclass{amsart}
\usepackage{amssymb} 
\usepackage{color}

 \newtheorem{theorem}{Theorem}[section]

 \newtheorem{proposition}[theorem]{Proposition}
 \newtheorem{lemma}[theorem]{Lemma}
 \newtheorem{corollary}[theorem]{Corollary}
 \newtheorem{remark}[theorem]{Remark}

 \numberwithin{equation}{section}
\def\reptoss{{\rm Rep}_2(m)_{ss}(\mathbb{C})}

\def\repss{{\rm Rep}_n(m)_{ss}(\mathbb{C})}
\def\repssk{{\rm Rep}_n(m)_{ss}(K)}
\def\repair{{\rm Rep}_n(m)_{air}(\mathbb{C})}
\def\reptoairab{{\rm Rep}_2(m)_{air}}
\def\repairk{{\rm Rep}_n(m)_{air}(K)}

\def\reptoair{{\rm Rep}_2(m)_{air}(\mathbb{C})}
\def\chair{{\rm Ch}_n(m)_{air}(\mathbb{C})}

\def\chtoair{{\rm Ch}_2(m)_{air}(\mathbb{C})}
\def\chtoairab{{\rm Ch}_2(m)_{air}}
\def\pglto{{\rm PGL}_2(\mathbb{C})}
\def\pgltok{{\rm PGL}_2(K)}
\def\pgl{{\rm PGL}_n(\mathbb{C})}
\def\pglk{{\rm PGL}_n(K)}
\def\configto{F_2(\mathbb{C}^m)}
\def\config{F_n(\mathbb{C}^m)}
\def\configk{F_n(K^m)}
\def\unconfig{\mathcal{C}_n(\mathbb{C}^m)}
\def\unconfigk{\mathcal{C}_n(K^m)}
\def\repunto{{\rm Rep}_2(m)_{u}(\mathbb{C})}
\def\repuntok{{\rm Rep}_2(m)_{u}(K)}
\def\chss{{\rm Ch}_n(m)_{ss}(\mathbb{C})}
\def\chssk{{\rm Ch}_n(m)_{ss}(K)}
\def\chtoss{{\rm Ch}_2(m)_{ss}(\mathbb{C})}

\def\chunto{{\rm Ch}_2(m)_{u}(\mathbb{C})}
\def\chuntok{{\rm Ch}_2(m)_{u}(K)}
\def\ssnm{{\rm SS}_n(m)}
\def\dn{{\rm D}_n}
\def\utom{{\rm U}_2(m)}
\def\nto{{\rm N}_2}
\def\bto{{\rm B}_2}
\def\tto{{\rm T}}
\def\tton{{\rm T}_n}
\def\ttotwo{{\rm T}_2}
\def\hn{{\rm H}_n}
\def\hto{{\rm H}_2}
\def\flag{{\rm Flag}(\mathbb{C}^n)}
\def\fq{\mathbb{F}_q}

\def\fqbar{\overline{\mathbb{F}}_q}
\def\qh{\displaystyle q^{\frac{1}{2}}}

\begin{document}

\title
[VHP of the moduli spaces of representations
     of degree $2$ for free monoids]
{
Virtual Hodge polynomials
of the moduli spaces of representations of degree $2$
for free monoids}
\author{Kazunori Nakamoto and Takeshi Torii}
\address{%Center for Life Science Research,
%Interdisciplinary Graduate School of Medicine and Engineering,
Center for Medical Education and Sciences,
Faculty of Medicine, 
University of Yamanashi,
Yamanashi 409--3898, Japan}
\email{nakamoto@yamanashi.ac.jp}

\address{Department of Mathematics, 
%Faculty of Sience, 
Okayama University,
Okayama 700--8530, Japan}
\email{torii@math.okayama-u.ac.jp}
\thanks{The first author was partially supported by 
JSPS KAKENHI Grant Number 23540044, 15K04814.
The second author was partially supported by
JSPS KAKENHI Grant Number 22540087, 25400092.}



\date{January 20, 2016\ ({\tt version~2.0})}

\begin{abstract}
In this paper we study the topology of the moduli spaces of representations
of degree $2$ for free monoids.
%We calculate the 
%(integral and rational)
%cohomology groups 
%of these moduli spaces.
We calculate the virtual Hodge 
polynomials of 
the character varieties  
for several types of $2$-dimensional representations.
%of degree $2$.
Furthermore,
we count the number of isomorphism classes for each type of 
$2$-dimensional representations 
%of degree $2$ 
over any finite field $\fq$,
and show that the number coincides with
the virtual Hodge polynomial evaluated at $q$. 
\end{abstract}

\maketitle

\section{Introduction}

The moduli spaces of representations have been studied 
in various contexts (see, for example,
\cite{FrickeKlein,King,N-S1,N-S2}).
In order to study the representation theory over schemes,
the first author defined the representation variety and
the character variety over ${\mathbb Z}$ of absolutely
irreducible representations 
for groups or monoids in \cite{Nakamoto1}.
Although they are of great importance in the representation
theory over schemes,
it is difficult to analyze them. 
To overcome this difficulty,
the first author 
introduced representations with Borel mold 
and studied the representation variety and the
character variety of representations with Borel mold
in \cite{Nakamoto2}.
Furthermore, we studied the topology of these moduli spaces 
of representations with Borel mold for free monoids
over the field $\mathbb{C}$ of complex numbers
in \cite{topos}.

Let us introduce a general framework for the representation theory
over schemes for free monoids (cf.~\cite{Nakamoto3,topos}).
Let ${\rm Rep}_n(m)=({\rm M}_n)^m$
the product of $m$-copies of the full $n\times n$ 
matrix ring,
which is the representation variety of degree $n$
over $\mathbb{Z}$
for the free monoid with $m$ generators. 
There is a decomposition of ${\rm Rep}_n(m)$
by locally closed subschemes
\[ {\rm Rep}_n(m) = \bigcup_{h=1}^{n^2} {\rm Rep}_n(m)_{{\rm rk}h},\]
where ${\rm Rep}_n(m)_{{\rm rk }h}$ is
the moduli space of representations which
generate a subalgebra of rank $h$ in the matrix algebra
${\rm M}_n$. Here we use ``rank'' in the sense of the rank of a 
subalgebra of the full matrix ring as projective module or 
locally free sheaf, because 
we deal with subalgebras (which is a projective module or locally free sheaf) 
of the matrix algebra ${\rm M}_n$ over 
an arbitrary commutative ring or scheme.   
The open subscheme ${\rm Rep}_n(m)_{{\rm rk}n^2}$ is the moduli
space ${\rm Rep}_n(m)_{air}$
of absolutely irreducible representations
\[ {\rm Rep}_n(m)_{air}={\rm Rep}_n(m)_{{\rm rk}n^2}.\]
In degree $2$ case we can study these moduli spaces
in more detail as in \cite{Nakamoto3}.
In particular, ${\rm Rep}_2(m)_{{\rm rk}3}$ is the moduli space
${\rm Rep}_2(m)_B$ of representations with Borel mold. 
We have two subschemes of ${\rm Rep}_2(m)_{{\rm rk}2}$.
One is the moduli space ${\rm Rep}_2(m)_{ss}$ of semi-simple representations,
and the other is the moduli space ${\rm Rep}_2(m)_{u}$ of
unipotent representations.
Then ${\rm Rep}_2(m)_u$ is closed in ${\rm Rep}_2(m)_{{\rm rk}2}$,
and ${\rm Rep}_2(m)_{ss}$ is its open complement.
We can regard ${\rm Rep}_2(m)_{{\rm rk}1}$ as
the moduli space ${\rm Rep}_2(m)_{sc}$
of scalar representations.
The conjugate action of matrices
induces an action of the group scheme 
${\rm PGL}_2$ on ${\rm Rep}_2(m)_{\ast}$ 
for $\ast = air, B, ss, u, sc$. 
The character variety ${\rm Ch}_2(m)_{\ast}$
is defined as
\[ {\rm Ch}_2(m)_{\ast}={\rm Rep}_2(m)_{\ast}/{\rm PGL}_2. \] 
%for $\ast = air, B, ss, u, sc$. 
In \cite{Nakamoto1}, \cite{Nakamoto2}, and \cite{Nakamoto3}, 
the first author showed that
${\rm Ch}_2(m)_{\ast}$ is a universal geometric quotient
of ${\rm Rep}_2(m)_{\ast}$ by ${\rm PGL}_2$ for $\ast = air, B, ss, u, sc$.  
(For $\ast = u$, we need to divide ${\rm Rep}_2(m)_{u}$ into two parts: 
the ${\Bbb Z}[1/2]$-part and the ${\Bbb F}_2$-part. 
More precisely, see \cite{Nakamoto3}.) 

In this paper we study the topology of 
the moduli spaces of representations of degree $2$ 
for free monoids.
We calculate the integral and rational cohomology groups 
of representation varieties ${\rm Rep}_2(m)_{\ast}$ and 
character varieties ${\rm Ch}_2(m)_{\ast}$ 
over the field $\mathbb{C}$ of complex numbers for $\ast = ss, u, sc$.  
Also, we give the virtual Hodge polynomial
of ${\rm Rep}_2(m)_{\ast}$ and ${\rm Ch}_2(m)_{\ast}$ 
for $\ast = ss, u, sc, air$. 
See \S\ref{section:VHP-air} 
for the definition of virtual Hodge polynomials.
%The following is one of the main theorems of this paper.
In this paper we obtain the following theorems:

\begin{theorem}[Theorem~\ref{thm:VHPreptoair}]
\label{theorem:Main-Theorem1}
The virtual Hodge polynomial of 
$\reptoair$
is given by
\[  \begin{array}{rcl}
     {\rm VHP}(\reptoair)(z)&=&
     (1-z^m)(1-z^{m-1}),\\[3mm]
     {\rm VHP}_c(\reptoair)(z)&=&
     z^{2m+1}(z^m-1)(z^{m-1}-1).\\
    \end{array} \]
\end{theorem}

%%{\color{red}
%The conjugate action of matrices
%induces an action of the group scheme 
%${\rm PGL}_n$ on ${\rm Rep}_n(m)_{air}$.
%The character variety ${\rm Ch}_n(m)_{air}$
%is defined as
%\[ {\rm Ch}_n(m)_{air}={\rm Rep}_n(m)_{air}/{\rm PGL}_n.\]
%In \cite{Nakamoto1} the first author showed that
%${\rm Ch}_n(m)_{air}$ is a universal geometric quotient
%of ${\rm Rep}_n(m)_{air}$ by ${\rm PGL}_n$.

%On the virtual Hodge polynomial of the character variety
%of absolutely irreducible representations of degree $2$,
%we have the following theorem.
%%}

%\begin{conjecture}[Conjecture~\ref{conj:VHPchtoair}]
\begin{theorem}[Theorem~\ref{thm:VHPchtoair}]
%\label{conjecture:Main-Theorem2}
\label{thm:Main-Theorem2}
The virtual Hodge polynomial of 
$\chtoair$ is given by 
\[ \begin{array}{rcl}
    {\rm VHP}(\chtoair)(z)&=&
   \frac{(1-z^m)(1-z^{m-1})}{1-z^2},\\[3mm]
    {\rm VHP}_c(\chtoair)(z)&=&
   \frac{z^{2m}(z^m-1)(z^{m-1}-1)}{z^2-1}.\\
   \end{array}\]
%\end{conjecture}
\end{theorem}
\noindent 
These results can be obtained by considering the mixed Hodge structures on 
the rational cohomology of ${\rm Rep}_2(m)_{\ast}$ 
for $\ast = B, ss, u, sc$ and 
calculating these virtual Hodge polynomials.

Furthermore,
we count the numbers of $\fq$-valued points of the
moduli spaces of representations of degree $2$
for free monoids.
In particular,
the number of $\fq$-valued points of
${\rm Rep}_2(m)_{air}$ is the number
of absolutely irreducible representations
of degree $2$ over $\fq$ for the free monoid
with $m$ generators,
and the number of $\fq$-valued points
of ${\rm Ch}_2(m)_{air}$ is the number
of isomorphism classes of such representations.

\begin{theorem}[cf.~Theorems~\ref{thm:number-fq-air} 
and \ref{th:numberofchair}]
\label{main-theorem-fq-counting}
For any finite field $\fq$,
the number of absolutely irreducible 
representations of degree $2$ over $\fq$
for the free monoid with $m$ generators 
is given by
\[ |{\rm Rep}_2(m)_{air}(\fq)|=
     {\rm VHP}_c({\rm Rep}_2(m)_{air}(\mathbb{C}))(q),\]
and the number of isomorphism classes of
such representations is given by
\[     |{\rm Ch}_2(m)_{air}(\fq)|=
     {\rm VHP}_c({\rm Ch}_2(m)_{air}(\mathbb{C}))(q). \]
\end{theorem}

Note that if there exists a polynomial with integer coefficients
which counts the number of $\fq$-valued points
of a separated scheme $X$ of finite type over $\mathbb{Z}$,
then $|X(\fq)|={\rm VHP}_c(X(\mathbb{C}))(q)$.
For more details,
see \cite[\S6]{HRV}. 
We should mention that these results coincide with the 
calculation of \cite[Example~7.2]{Reineke}. 
Reineke has already calculated the virtual Hodge 
polynomials of ${\rm Ch}_2(m)_{air}(\mathbb{C})$.

By these results, we have the following main theorem of this paper:  

\begin{theorem}[Theorem~\ref{thm:VHPch}]\label{theorem:Main-Theorem4} 
The number of ${\Bbb F}_q$-valued points of ${\rm Ch}_2(m)$ is given by 
\[ 
\begin{array}{ccl} 
|{\rm Ch}_2(m)({\mathbb F}_q)| & = & |{\rm Ch}_2(m)_{\rm air}({\mathbb F}_q)| 
+ |{\rm Ch}_2(m)_{ss}({\mathbb F}_q)| 
+ |{\rm Ch}_2(m)_{sc}({\mathbb F}_q)| \vspace*{1ex} \\ 
 & = & \displaystyle \frac{q^{2m+2}(q^{2m-3} - q^{m-2} - q^{m-3} + 1)}{q^2 - 1}.   
\end{array} 
\] 
In particular, the virtual Hodge polynomial of 
${\rm Ch}_2(m)$ is given by 
\[ 
%\begin{array}{ccl} 
{\rm VHP}_{c}({\rm Ch}_2(m))(z) 
  =  \displaystyle \frac{z^{2m+2}(z^{2m-3} - z^{m-2} - z^{m-3} + 1)}{z^2 - 1}.   
%\end{array} 
\]
\end{theorem}

In the present paper, we deal with only $2$-dimensional representations of 
free monoids. However, our strategy is available for other groups or monoids. 
We can calculate the numbers of equivalence classes of 
absolutely irreducible $2$-dimensional representations 
over ${\Bbb F}_q$ for finitely generated groups or monoids $\Gamma$ by 
calculating the number of the ${\Bbb F}_q$-valued points of 
the representation variety ${\rm Rep}_2(\Gamma)$ and 
the others ${\rm Rep}_2(\Gamma)_{\ast}$ for 
$\ast = B, s.s., u, sc$. 
If there are polynomials with integer coefficients which count 
the numbers of the ${\Bbb F}_q$-valued points of the reprensetation variety 
${\rm Rep}_2(\Gamma)$ and the others ${\rm Rep}_2(\Gamma)_{\ast}$, then 
the polynomials coincide with virtual Hodge polynomials 
of the corresponding schemes. 

For $3$-dimensional representations, we can also calculate 
the virtual Hodge polynomials of the 
character varieties of absolutely irreducible representations 
for finitely generated free monoids. Indeed, there are $26$ types of 
subalgebras of the full matrix ring of degree $3$, and we can 
calculate the virtual Hodge polynomial of the representation variety  
${\rm Rep}_3(m)_{\ast}$ associated to each subalgebra of ${\rm M}_3$  
in the same way as $2$-dimensional representations. 

Here we should point out the results of \cite{SL2-ch}, \cite{SL3-ch}, 
\cite{Moz-number}, \cite{Moz-arith}, \cite{Reineke},   
and so on. After we wrote the present paper, 
we found out their papers. Their results are very much related to our paper. 
Our strategy to calculate the numbers of ${\Bbb F}_q$-valued points 
and the virtual Hodge polynomials of the 
representation varieties and the character varieties is essentially 
same as \cite{SL2-ch} and \cite{SL3-ch}.  
Moreover, it is much easier to calculate the virtual Hodge 
polynomials of ${\rm Ch}_n(m)_{air} := {\rm Rep}_n(m)_{air}/{\rm PGL}_n$ by 
the method of \cite{Reineke} than 
to calculate them by our strategy. 
However, we believe that it is worth publishing our results. 
Not only are our objects ${\rm Rep}_2(m)_{\ast}$ and ${\rm Ch}_2(m)_{\ast}$ 
different from the ${\rm SL}_2({\Bbb C})$-character 
varieties of free groups, but also we 
give a geometric meaning to their stratifications of 
the representation varieties in \cite{SL2-ch}.   
Each stratification represents a certain moduli functor 
which is described in terms of representation theory.  
This interpretation allows us to overview the representation varieties 
and the character varieties from viewpoints of 
algebraic geometry, algebraic topology, representation theory, and so on.  
Furthermore, as far as we know, our strategy is 
the only way to calculate the virtual Hodge polynomial of 
${\rm Ch}_2(m)$.

The organization of this paper is as follows:
In \S\ref{section:ss-representaion} 
we study the representation variety and the
character variety of semi-simple representations.
We give descriptions of these moduli spaces.
Then we calculate the integral and rational cohomology
groups of them.
In \S\ref{section:unipotent-representations} 
we define the representation variety and
the character variety of unipotent representations
of degree $2$.
Then we give descriptions and calculate
the cohomology groups of them.
In \S\ref{section:scalar-representations}
we define and study
the moduli spaces of scalar representations.
In \S\ref{section:VHP-air}
we study the virtual Hodge polynomials
of the representation varieties
and the character varieties of representations of degree $2$
for free monoid.
Then we prove Theorem~\ref{theorem:Main-Theorem1}
and Theorem~\ref{thm:Main-Theorem2}.
%and give Conjecture~\ref{conjecture:Main-Theorem2}.
In \S\ref{section:counting}
we count the number of $\fq$-valued points
of the moduli spaces,
and prove Theorem~\ref{main-theorem-fq-counting}
 and Theorem~\ref{theorem:Main-Theorem4}. 

Let $\mathbb{Z}$, $\mathbb{Q}$,
$\mathbb{R}$, $\mathbb{C}$
be the ring of integers,
the field of rational numbers,
the field of real numbers,
the field of complex numbers,
respectively.
In this paper we denote by
$H^*(X)$ the integral cohomology groups,
and by $H^*(X;\mathbb{Q})$
the rational cohomology groups of a space $X$.
For a graded module $V$ over $\mathbb{Z}$ or
$\mathbb{Q}$,
we denote by $\Lambda(V)$ the free commutative graded algebra
on $V$.

\section{The moduli spaces of semi-simple representations}
\label{section:ss-representaion}

In this section we study the moduli spaces related to semi-simple
representations. We give descriptions for the moduli spaces
and calculate the cohomology groups of them.

\subsection{Descriptions for the moduli spaces
of semi-simple representations}

Let $K$ be an algebraically closed field.
For $(A_1,\ldots,A_m)\in ({\rm M}_n(K))^m$,
we write $A_k=(a_{ij}(k))_{ij}$.
Then we define $a_{ij}\in K^m$
by $a_{ij}=(a_{ij}(1),\ldots,a_{ij}(m))$.
We denote by $\dn$ the $K$-subalgebra 
of ${\rm M}_n(K)$
consisting of diagonal matrices.
Let $\ssnm$ be the subspace of $(\dn)^m$ given by
\[ \ssnm =\{(A_1,\ldots,A_m)\in (\dn)^m|\
           \mbox{\rm $A_1,\ldots,A_m$ generate $\dn$ as a 
                     $K$-algebra} \}. \]
We define the moduli space $\repssk$ of semi-simple
representations by
\[ \begin{array}{cl}
      &\repssk\\[3mm] 
    =&\{(A_1,\ldots,A_m)\in ({\rm M}_n(K))^m|\
           \mbox{\rm $(PA_1P^{-1},\ldots,PA_mP^{-1})\in \ssnm$ for
             some $P\in \pglk$}\}. 
  \end{array}\]
Note that $\ssnm$ is a subspace of $\repssk$.
Let $\configk$ be the configuration
space of distinct ordered  $n$-points in $K^m$:
\[ \configk = \{(p_1,\ldots,p_n)\in (K^m)^n|\
    \mbox{\rm $p_i\neq p_j$ for $i\neq j$}\}.  \]
For $p=(p_1,\ldots,p_n)\in (K^m)^n$,
we set 
\[ \varphi(p)= (A_1,\ldots,A_m)\in (\dn)^m\]
where $a_{ij}= p_i$ if $i=j$, and
      $a_{ij}=0$    if $i\neq j$.
Then we obtain an isomorphism  
$\varphi: (K^m)^n\stackrel{\cong}{\longrightarrow} (\dn)^m$.
By
\cite[Lemma~3.3]{topos},
we obtain the following lemma.

\begin{lemma}\label{lemma:config-ssnm}
The map $\varphi$ induces an isomorphism
$\configk\stackrel{\cong}{\longrightarrow} \ssnm$
of smooth algebraic varieties.
So we can regard $\configk$ as a subspace of 
$\repssk$.
\end{lemma}

Let $\Sigma_n$ be the symmetric group on $n$-letters.
We regard $\Sigma_n$ as a subgroup of $\pglk$
consisting of permutation matrices.
Let $\tton$ be the diagonal subgroup of $\pglk$.
We denote by $\hn$ the subgroup of $\pglk$ 
generated by $\Sigma_n$ and $\tton$.
Then $\tton$ is a normal subgroup of $\hn$
and $\hn$ is isomorphic to the semi-direct product 
$\tton\rtimes\Sigma_n$.
It is easy to prove the following lemma.

\begin{lemma}\label{lemma:normalizer-ss}
If $P\in\pglk$ satisfies $P\dn P^{-1}=\dn$,
then $P\in \hn$.
\end{lemma}

Since $\pglk$ acts on 
$\repssk$ by conjugation,
we can extend the inclusion map $\varphi: \configk\to 
\repssk$
to a map $\pglk\times 
\configk\longrightarrow \repssk$.
Note that the subgroup $\hn$ preserves the subspace 
$\configk$.
So this map factors through $\pglk
\times_{\hn} \configk$.

\begin{theorem}\label{theorem:repss-description}
There is an isomorphism of smooth algebraic varieties 
with $\pglk$-action
\[ \repssk\cong \pglk\times_{\hn}\configk.\] 
In particular, 
$\dim_K\repssk = mn+n(n-1)$.
\end{theorem}

\proof
For any $A=(A_1,\ldots,A_m)\in\repssk$,
there exists $P\in\pglk$ such that
$PAP^{-1}\in \configk$.
So $\pglk\times_{\hn} 
\configk\to \repssk$ is surjective.
Suppose that $PUP^{-1}=QVQ^{-1}$ for
$P,Q\in\pglk$, $U,V\in \configk$.
Then $U=P^{-1}QVQ^{-1}P\in \configk$.
So both of $U$ and $P^{-1}QVQ^{-1}P$
generate $\dn$. 
By Lemma~\ref{lemma:normalizer-ss},
$P^{-1}Q=B\in \hn$.
Then we see that
$(Q,V)=(PB,V)=(P,BVB^{-1})=(P,U)$
in $\pglk\times_{\hn}\configk$.
This completes the proof.
\qed

\

\begin{remark}\rm
There exists a scheme ${\rm Rep}_2(m)_{ss}$
of finite type over $\mathbb{Z}$ by \cite{Nakamoto3},
and ${\rm Rep}_2(m)_{ss}(K)$ is the associated algebraic variety.
\end{remark}

By Theorem~\ref{theorem:repss-description},
we have $\repssk\cong \pglk/
\tton\times_{\Sigma_n}\configk$.
Note that $\pglk/\tton$
is the space of $n$-tuples $(l_1,\ldots,l_n)$,
where $l_i$ is a $1$-dimensional subspace of $K^n$
for $i=1,\ldots,n$,
and $\sum_{i=1}^n l_i=K^n$.
Then $\Sigma_n$ acts on $\pglk/\tton$ as permutations
of lines.

We define the character variety $\chssk$ of
semi-simple representations as 
the quotient space of $\repssk$ 
by $\pglk$: 
\[ \chssk = \repssk/\pglk. \]
Let $\unconfigk$ be the configuration space of 
distinct unordered $n$-points in $K^m$
which is defined to be the quotient space
$\configk/\Sigma_n$:
\[ \unconfigk= \configk/\Sigma_n.\]
Since $\chssk = \repssk/\pglk\cong 
\configk/\hn$
and the subgroup $\tton$ acts trivially on $\configk$,
we obtain the following corollary.

\begin{corollary}\label{cor:chss-unconfig}
There is an isomorphism of smooth algebraic varieties
\[ \chssk\cong \unconfigk. \]
\end{corollary}

\begin{remark}\rm
There exists a scheme ${\rm Ch}_2(m)_{ss}$
of finite type over $\mathbb{Z}$ by \cite{Nakamoto3},
and ${\rm Ch}_2(m)_{ss}(K)$ is the associated algebraic variety.
\end{remark}

\subsection{Integral cohomology groups of the moduli spaces of degree $2$}
\label{subsection:ssintegralcohomologyofdegree2}

In this subsection we restrict our attention to
the degree $2$ representations
over $\mathbb{C}$.
We study the integral cohomology 
groups of the moduli spaces related to semi-simple representations
of degree $2$.

We denote by $S^r$ 
the $r$-sphere,
and by ${\rm P}^r$
the projective $r$-space. 

\begin{lemma}\label{lemma:config-homotopy-type}
The configuration space $\configto$
with $\Sigma_2$-action
is equivariantly homotopy equivalent 
to $S^{2m-1}$ with antipodal $\Sigma_2$-action.  
Hence $\mathcal{C}_2(\mathbb{C}^m)$ is homotopy equivalent
to ${\rm P}^{2m-1}(\mathbb{R})$.
\end{lemma}

\proof
We regard $S^{2m-1}$ as the unit sphere in $\mathbb{C}^m$.
We have a $\Sigma_2$-equivariant map
from $\configto$ to $S^{2m-1}$ given by
$(p_1,p_2)\mapsto \frac{p_1-p_2}{||p_1-p_2||}$,
where $||\cdot||$ is the standard norm in $\mathbb{C}^m$.
It is easy to see that this map is a non-equivariant
homotopy equivalence.
Since $\configto$ and $S^{2m-1}$ are free $\Sigma_2$-spaces,
this map is a $\Sigma_2$-equivariant homotopy equivalence.
\qed

\

By Corollary~\ref{cor:chss-unconfig} and 
Lemma~\ref{lemma:config-homotopy-type},
we can calculate the cohomology groups
of $\chtoss$.

\begin{corollary}
We have an isomorphism of commutative graded  algebras
\[ H^*(\chtoss)\cong \Lambda(s,t)/(2s,s^m,st)\]
with $|s|=2$ and $|t|=2m-1$.
\end{corollary}

\begin{lemma}\label{lemma:flag-homotopy-type}
The $\Sigma_2$-space $\pglto/\ttotwo$ is equivariantly
homotopy equivalent to 
$S^2$ with antipodal $\Sigma_2$-action.
Hence $\pglto/\hto$ is homotopy equivalent to
${\rm P}^2(\mathbb{R})$.
\end{lemma}

\proof
Let $X$ be the space of pairs $(L_1,L_2)$
where $L_1$ and $L_2$ are orthogonal lines in $\mathbb{C}^2$.
Then $X$ is a $\Sigma_2$-subspace of $\pglto/\ttotwo$,
and the inclusion $X\hookrightarrow \pglto/\ttotwo$ is non-equivariantly
homotopy equivalent.
Since $\pglto/\ttotwo$ and $X$ are free $\Sigma_2$-spaces,
$X$ is $\Sigma_2$-equivariantly equivalent to $\pglto/\ttotwo$.

We regard ${\rm P}^1(\mathbb{C})$ as the space
of lines in $\mathbb{C}^2$.
Since the orthogonal complement $L^{\perp}$
of $L\in {\rm P}^1(\mathbb{C})$
is uniquely determined,
we can identify $X$ with ${\rm P}^1(\mathbb{C})$.
Then the nontrivial element $\tau$ in $\Sigma_2$ acts
on ${\rm P}^1(\mathbb{C})$ as $\tau(L)=L^{\perp}$
for $L\in {\rm P}^1(\mathbb{C})$.
The space ${\rm P}^1(\mathbb{C})$ with this $\Sigma_2$-action
can be identified with $S^2$ with antipodal $\Sigma_2$-action.
\qed

\begin{corollary}\label{cor:reptoss-description}
The space $\reptoss$ is homotopy equivalent to
$S^2\times_{\Sigma_2}S^{2m-1}$.
\end{corollary}

\proof
By Theorem~\ref{theorem:repss-description},
$\reptoss\cong \pglto/\ttotwo\times_{\Sigma_2}\configto$.
Then the corollary follows from
Lemmas~\ref{lemma:config-homotopy-type} and 
\ref{lemma:flag-homotopy-type}.
\qed

\begin{proposition}\label{prop:coh-rep_m=1}
We have an isomorphism of commutative graded  algebras
\[ H^*({\rm Rep}_2(1)_{ss}(\mathbb{C}))\cong \Lambda (a,b)/(2b,ab)\]
%\[ H^*(S^2\times_{\Sigma_2}S^1)\cong \wedge(a,b)/(2b,ab)\]
with
$|a|=1$ and $|b|=3$.   
\end{proposition}

\proof
By Corollary~\ref{cor:reptoss-description},
${\rm Rep}_2(1)_{ss}(\mathbb{C})$ is homotopy equivalent to
$S^2\times_{\Sigma_2}S^1$.
Consider the Serre spectral sequence
associated with the fibration
$S^2\to S^2\times_{\Sigma_2}S^1\to S^1/\Sigma_2\cong S^1$:
\[ E_2^{p,q}=H^p(\mathbb{Z};H^q(S^2))\Longrightarrow 
   H^{p+q}(S^2\times_{\Sigma_2}S^1).\]
Note that $\pi_1(S^1)\cong\mathbb{Z}$ acts 
nontrivially on $H^2(S^2)$.
Then we have
$E_2^{0,0}\cong E_2^{1,0}\cong \mathbb{Z}$,
$E_2^{1,2}\cong \mathbb{Z}/2$,
and $E_2^{p,q}=0$ otherwise. 
Hence the spectral sequence collapses.
Since there are no extension problems,
we obtain the proposition.
\qed

\begin{theorem}\label{thm:coh-rep_m}
For $m\ge 2$, we have
an isomorphism of commutative graded  algebras
\[ H^*(\reptoss)\cong
   \Lambda (a,b)/(2a,a^2)\]
%\[ H^*(S^2\times_{\Sigma_2}S^{2m-1})\cong
%   \wedge (a,b)/(2a),\]
with
$|a|=2$ and $|b|=2m-1$. 
\end{theorem}

\if0
\proof
By Corollary~\ref{cor:reptoss-description},
$\reptoss$ is homotopy equivalent to $S^2\times_{\Sigma_2}S^{2m-1}$.
Let $B\Sigma_2$ be the classifying space of $\Sigma_2$.
Consider the Serre spectral sequence
associated with the fibration
$S^2\times S^{2m-1}\to S^2\times_{\Sigma_2}S^{2m-1}\to B\Sigma_2$:
\[ E_2^{p,q}=H^p(\Sigma_2;H^q(S^2\times S^{2m-1}))\Longrightarrow
             H^{p+q}(S^2\times_{\Sigma_2}S^{2m-1}).\]
Note that $\Sigma_2$ acts nontrivially on $H^2(S^2)$,
and trivially on $H^{2m-1}(S^{2m-1})$.
Then we have $E_2^{0,0}\cong E_2^{0,2m-1}\cong \mathbb{Z}$,
$E_2^{2n,0}\cong E_2^{2n-1,2}\cong E_2^{2n,2m-1}\cong
E_2^{2n-1,2m+1}\cong \mathbb{Z}/2$ for $n>0$,
and $E_2^{p,q}=0$ otherwise.
By degree reasons, we see that $d_2=0$. 
By comparing the spectral sequence associated 
with the fibration $S^2\to S^2/\Sigma_2\cong \mathbb{RP}^2\to
B\Sigma_2$, we obtain that $d_3: E_3^{2n-1,2}\to E_3^{2n+2,0}$
and $d_3: E_3^{2n-1,2m+1}\to E_3^{2n+2,2m-1}$
are isomorphisms for $n>0$.
Hence $E_4^{0,0}\cong E_4^{0,2m-1}\cong \mathbb{Z}$,
$E_4^{2,0}\cong E_4^{2,2m-1}\cong \mathbb{Z}/2$,
and $E_4^{p,q}=0$ otherwise.
By degree reasons,
the spectral sequence collapses from the $E_4$-term.
Since there are no extension problems,
we obtain the theorem.
\qed
\fi

\proof
By Corollary~\ref{cor:reptoss-description},
$\reptoss$ is homotopy equivalent to $S^2\times_{\Sigma_2}S^{2m-1}$.
Consider the Serre spectral sequence
associated with the fibration
$S^{2m-1}\to S^2\times_{\Sigma_2}S^{2m-1}\to S^2/\Sigma_2\cong 
{\rm P}^2(\mathbb{R})$:
\[ E_2^{p,q}=H^p({\rm P}^2(\mathbb{R});\mathcal{H}^q(S^{2m-1}))\Longrightarrow
             H^{p+q}(S^2\times_{\Sigma_2}S^{2m-1}),\]
where $\mathcal{H}^q(S^{2m-1})$ is the local coefficient system
determined by the action of the 
fundamental group $\pi_1({\rm P}^2(\mathbb{R}))$
on $H^q(S^{2m-1})$.
In this case $\pi_1({\rm P}^2(\mathbb{R}))\cong\mathbb{Z}/2$
acts trivially on $H^{2m-1}(S^{2m-1})$.
Then we have $E_2^{0,0}\cong E_2^{0,2m-1}\cong \mathbb{Z}$,
$E_2^{2,0}\cong E_2^{2,2m-1}\cong\mathbb{Z}/2$,
and $E_2^{p,q}=0$ otherwise.
Hence the spectral sequence collapses.
Since there are no extension problems,
we obtain the theorem.
\qed

\if0
Hence we can interpret that
$\repss$ is the space of unordered
$n$-tuples of pairs $\{(p_1,l_1),\ldots,(p_n,l_n)\}$
such that
$p_i\in\mathbb{C}^n\ (i=1,\ldots,n),
p_i\neq p_j\ (i\neq j)$,
and 
$l_i$ is a line in $\mathbb{C}^n$
such that 
$l_1+\cdots +l_n=\mathbb{C}^n$.

For a space $X$,
we let $\mathcal{F}(X)$ be the set of 
all $n$-tuples $(L_1,\ldots,L_n)$,
where $L_i$ is a complex line bundle 
over $X$ for $i=1,\ldots,n$,
and $L_1\oplus\cdots\oplus L_n$
is trivial.
Then $\mathcal{F}$ is a functor
from the homotopy category of (nice) spaces
to the category of sets.

\begin{lemma}
The functor $\mathcal{F}$ is represented by $\pgl/\tton$.
\end{lemma}
\fi

%\input{rational}
\subsection{Rational cohomology groups of the moduli spaces}

In this subsection we study the rational cohomology groups
of the moduli spaces related to semi-simple representations
over $\mathbb{C}$.
%We denote by $H^*(X;\mathbb{Q})$ the cohomology group
%of a space $X$ with coefficients in the rational number field $\mathbb{Q}$.

Recall that $\pgl/\tton$ is the space of ordered
$n$-lines $(l_1,\ldots,l_n)$ in $\mathbb{C}^n$
such that $\sum_{i=1}^nl_i=\mathbb{C}^n$.
We let $F_j=\sum_{i=1}^jl_j$ be the subspace of $\mathbb{C}^n$
spanned by $l_i$ for $1\le i\le j$.
Then $(F_1,\ldots,F_n)$ is a complete flag in $\mathbb{C}^n$.
We denote by $\flag$ the flag variety,
which is the space of complete flags in $\mathbb{C}^n$. 
So we obtain a map
$\pgl/\tton\to\flag$ of complex manifolds.
Since this map is a homotopy equivalence,
we obtain the following lemma.

\begin{lemma}
The cohomology group of $\pgl/\tton$ is given by
\[ H^*(\pgl/\tton)\cong \mathbb{Z}[t_1,\ldots,t_n]/(c_1,\ldots,c_n),\]
where
$|t_1|=\cdots =|t_n|=2$, 
and $c_i$ is the $i$th elementary symmetric polynomial
of $t_1,\ldots, t_n$ for $i=1,\ldots,n$.
The action of $\Sigma_n$ on $\pgl/\tton$ induces
an action on $H^*(\pgl/\tton)$, which is given by
permutations of $t_1,\ldots,t_n$.
\end{lemma}

\begin{lemma}\label{lemma:flag-regular-rep}
The rational cohomology $H^*(\pgl/\tton;\mathbb{Q})$
is the regular representation of $\Sigma_n$.
\end{lemma}

\proof
Since $\Sigma_n$ freely acts on $\pgl/\tton$,
there are no fixed points.
Let $\chi$ be the character of the representation
defined by $H^*(\pgl/\tton;\mathbb{Q})$.
By the Lefschetz fixed point formula,
$\chi(g)=0$ if $g$ is not the identity in $\Sigma_n$,
and $\chi(g)=n!$ if $g$ is the identity.
Hence $H^*(\pgl/\tton;\mathbb{Q})$ is the regular representation
of $\Sigma_n$.
\qed

%Recall that 
%$\repss\cong
%\pgl/{\tton}\times_{\Sigma_n}\config$
%by Theorem~\ref{theorem:repss-description}.

\

We can easily describe the rational cohomology groups of 
a quotient space by a free action of a finite group.
We put the following well-known lemma for 
the reader's convenience.

\begin{lemma}\label{lemma:rational-cohomology-finite-group-quotient}
Let $G$ be a finite group.
Suppose that $G$ freely acts on a space $X$ and that
the quotient map $X\to X/G$ is a principal $G$-bundle.
Then we have an isomorphism of commutative graded algebras
\[ H^*(X/G;\mathbb{Q})\cong H^*(X;\mathbb{Q})^G.\]
\end{lemma}

\proof
By the assumptions,
we have a fibration $X\to X/G\to BG$,
where $BG$ is the classifying space of $G$.
We consider the associated Serre spectral sequence
\[ E_2^{p,q}=H^p(G;H^q(X;\mathbb{Q}))
   \Longrightarrow H^{p+q}(X/G;\mathbb{Q}).\]
For any $G$-module $M$ over $\mathbb{Q}$,
we have $H^p(G;M)=0$ for $p>0$
(see, for example, \cite[{Chapter~III, Corollary~10.2}]{Brown}).
Hence the spectral sequence collapses
and we obtain that
$H^p(X/G;\mathbb{Q})\cong H^p(X;\mathbb{Q})^G$.
%This completes the proof.
\qed
 
\begin{theorem}\label{thm:rational-repnmss}
We have an isomorphism of commutative graded algebras
\[ H^*(\repss;\mathbb{Q})\cong 
  \left(H^*(\pgl/\tton;\mathbb{Q})\otimes
         H^*(\config;\mathbb{Q})\right)^{\Sigma_n}.\]
\end{theorem}

\proof
The symmetric group $\Sigma_n$ freely acts on
$\pgl/\tton\times \config$ and
the quotient space $\pglto/\tton\times_{\Sigma_n}\config$
is isomorphic to $\repss$
by Theorem~\ref{theorem:repss-description}.
Hence we obtain the theorem by
Lemma~\ref{lemma:rational-cohomology-finite-group-quotient}. 
%This implies that the rational cohomology 
%$H^*(\repss;\mathbb{Q})$ is isomorphic to
%the $\Sigma_n$-invaritant part of 
%$H^*(\pglto/\tton\times \config;\mathbb{Q})$.
\qed

\

Note that the cohomology groups of the configuration space
$\config$ is given by
\[ H^*(\config)=\Lambda(s(i,j))_{1\le i<j\le n}/I,\]
where $|s(i,j)|=2m-1$ for $1\le i< j\le n$ and
the ideal $I$ is generated by
\[   s(i,k)s(j,k)-s(i,j)s(j,k)+s(i,j)s(i,k) \]
for $1\le i<j<k\le n$.

\begin{corollary}\label{cor:rationaliso-config-repss}
For $n=2$,
the inclusion
$\configto\hookrightarrow \reptoss$ induces
an isomorphism of rational cohomology groups
\[ H^*(\reptoss;\mathbb{Q})\stackrel{\cong}{\longrightarrow}
   H^*(\configto;\mathbb{Q}),\]
which is an isomorphism of commutative graded algebras.
\end{corollary}

\proof
By Lemma~\ref{lemma:flag-homotopy-type},
we see that $\Sigma_2$ non-trivially 
acts on $H^2(\pglto/\tto_2;\mathbb{Q})$.
On the other hand,
$\Sigma_2$ trivially
acts on $H^*(\configto;\mathbb{Q})$
by Lemma~\ref{lemma:config-homotopy-type}.
These imply that
$(H^*(\pglto/\tto_2;\mathbb{Q})\otimes
H^*(\configto;\mathbb{Q}))^{\Sigma_2}
\cong H^*(\configto;\mathbb{Q})$.
\qed

\begin{remark}\rm
Corollary~\ref{cor:rationaliso-config-repss} 
also follows from 
Proposition~\ref{prop:coh-rep_m=1} and
Theorem~\ref{thm:coh-rep_m}.
%the results in 
%\S\ref{subsection:ssintegralcohomologyofdegree2}.
\end{remark}

\begin{corollary}
We have
\[ \dim H^{\rm even}(\repss;\mathbb{Q})=
   \dim H^{\rm odd}(\repss;\mathbb{Q})=\frac{n!}{2}.\] 
\end{corollary}

\proof
By induction on $n$,
we see that
$\dim H^{\rm even}(\config;\mathbb{Q})=\dim
H^{\rm odd}(\config;\mathbb{Q})= n!/2$.
Since $\mathbb{Q}[\Sigma_n]\otimes V\cong 
\mathbb{Q}[\Sigma_n]^{\oplus\dim V}$ as $\mathbb{Q}[\Sigma_n]$-modules
for any representation $V$ of $\Sigma_n$, 
we see that 
$H^*(\repss;\mathbb{Q})\cong H^*(\config;\mathbb{Q})$
as $\mathbb{Q}$-vector spaces for $*={\rm even},{\rm odd}$ 
by Lemma~\ref{lemma:flag-regular-rep} and 
Theorem~\ref{thm:rational-repnmss}.
\qed

\if0
\proof
By \cite[Theorem~5.1]{Cohen-Talyor},
$H^{\rm even}(\config;\mathbb{Q})$ and $H^{\rm odd}(\config;\mathbb{Q})$
are isomorphic to $1_{(1,2)}|^{\Sigma_n}$ as
$\mathbb{Q}[\Sigma_n]$-modules
respectively,
where $1_{(1,2)}$ is the trivial representation 
of the subgroup generated by the transposition $(1,2)$,
and $1_{(1,2)}|^{\Sigma_n}$ is the induced representation of 
$1_{(1,2)}$ to $\Sigma_n$.
Since $\mathbb{Q}[\Sigma_n]\otimes V\cong 
\mathbb{Q}[\Sigma_n]^{\oplus\dim V}$ as $\mathbb{Q}[\Sigma_n]$-modules
for any representation $V$ of $\Sigma_n$, 
we see that 
$H^*(\repss;\mathbb{Q})\cong H^*(\config;\mathbb{Q})$
as $\mathbb{Q}$-vector spaces for $*={\rm even},{\rm odd}$ 
by Lemma~\ref{lemma:flag-regular-rep} and 
Theorem~\ref{thm:rational-repnmss}.
\qed
\fi

\begin{lemma}\label{lemma:cohomology-unconfig}
The quotient map
$\config\to\unconfig$ 
induces an injection of rational cohomology
groups
\[ H^*(\unconfig;\mathbb{Q})\hookrightarrow 
   H^*(\config;\mathbb{Q}),\]
and the image of the map is identified with the $\Sigma_n$-invariant
submodule of $H^*(\config;\mathbb{Q})$.
We have an isomorphism of commutative graded algebras
\[ H^*(\unconfig;\mathbb{Q})\cong \Lambda(s)\]  
with
$|s|=2m-1$.
\end{lemma}

\proof
Since $\Sigma_n$ freely acts on $\config$,
we obtain that
$H^p(\unconfig;\mathbb{Q})\cong H^p(\config;\mathbb{Q})^{\Sigma_n}$
by Lemma~\ref{lemma:rational-cohomology-finite-group-quotient}.
Let $s=\sum_{i<j}s(i,j)\in H^{2m-1}(\config;\mathbb{Q})$.
Then we have
$H^*(\config;\mathbb{Q})^{\Sigma_n}=\Lambda(s)$
by \cite[Corollary~5.2]{Cohen-Talyor}.
This completes the proof.
\if
Since the action of $\Sigma_n$ on $\config$
is free, we have a fibre bundle
$\config\to\unconfig\to B\Sigma_n$.
We consider the associated Serre spectral sequence
\[ E_2^{p,q}=H^p(\Sigma_n;H^q(\config;\mathbb{Q}))
   \Longrightarrow H^{p+q}(\unconfig;\mathbb{Q}).\]
For any $\Sigma_n$-module over $\mathbb{Q}$,
we have $H^p(\Sigma_n;M)=0$ for $p>0$
(see, for example, \cite[{Chapter~III, Corollary~10.2}]{Brown}).
Hence the spectral sequence collapses
and we obtain that
$H^p(\unconfig;\mathbb{Q})\cong H^p(\config;\mathbb{Q})^{\Sigma_n}$.
Let $s=\sum_{i<j}s(i,j)\in H^{2m-1}(\config;\mathbb{Q})$.
Then we have
$H^*(\config;\mathbb{Q})^{\Sigma_n}=\Lambda(s)$
by \cite[Corollary~5.2]{Cohen-Talyor}.
This completes the proof.
\fi
\qed

\
 
By Corollary~\ref{cor:chss-unconfig},
we obtain the following proposition.

\begin{proposition}\label{prop:chss-rational}
The composition map $\config\to \repss\to \chss$
induces an injection of
the rational cohomology groups
\[ H^*(\chss;\mathbb{Q})\hookrightarrow H^*(\config;\mathbb{Q}).\]
%We denote by $s(i,j)$ the additive basis of $H^{2m-1}(\config)$,
%and we let $s=\sum_{i<j}s(i,j)$.
We have an isomorphism of commutative graded  algebras
\[ H^*(\chss;\mathbb{Q})\cong \Lambda (s)\]
with
$|s|=2m-1$.
\end{proposition}

\if
\proof
By Corollary~\ref{cor:chss-unconfig},
$\chss$ is isomorphic to
$\unconfig$ the configuration space of distinct 
unordered $n$-points in $\mathbb{C}^m$.
Hence 
$H^*(\chss;\mathbb{Q})$ is the $\Sigma_n$-invariant part
of $H^*(\config;\mathbb{Q})$.
Then the proposition follows from
\cite[Corollary~5.2]{Cohen-Talyor}.
\qed
\fi

\section{The moduli spaces of unipotent representations of degree $2$}
\label{section:unipotent-representations}

In this section we study the moduli spaces related to unipotent
representations of degree $2$. 
We give descriptions for the moduli spaces
and calculate the cohomology groups of them.
 
\subsection{Descriptions for the moduli spaces 
of unipotent representations of degree $2$}

Let $K$ be an algebraically closed field.
Let $\nto$ be the $K$-subalgebra of 
${\rm M}_2(K)$
generated by the following matrix
\[  \left(\begin{array}{cc}
           0 & 1\\
           0 & 0\\
    \end{array}\right).\]
Note that $\dim_{K}\nto=2$.
Let $\utom(K)$ 
be the subspace of $(\nto)^m$ given by
\[ \utom(K)
   =\left\{ (A_1,\ldots,A_m)\in (\nto)^m|\
    \mbox{\rm $A_1,\ldots,A_m$ generate $\nto$ as a $K$-algebra}
    \right\}.\]
Note that $\utom(K)$ is an algebraic variety associated
to a scheme $\utom$ over $\mathbb{Z}$. 
We define the moduli space $\repuntok$ of unipotent
representations of degree $2$ by
%to be
\[ \repuntok
    =\left\{(A_1,\ldots,A_m)\in ({\rm M}_2(K))^m\left|\
           \begin{array}{l}
           (PA_1P^{-1},\ldots,PA_mP^{-1})\in 
           \utom(K)\\[2mm] 
           \mbox{\rm for some $P\in\pgltok$}
           \end{array}\right.\right\}. \]
%\[ \begin{array}{rl}
%     &\repuntok\\[3mm] 
%    =&\{(A_1,\ldots,A_m)\in ({\rm M}_2(K))^m|\
%           \mbox{\rm $(PA_1P^{-1},\ldots,PA_mP^{-1})\in 
%           \utom(K)$
%           for some $P\in\pgltok$}\}. 
%  \end{array}\]
%\[ \{(A_1,\ldots,A_m)\in ({\rm M}_2(K))^m|\
%           \mbox{\rm $(PA_1P^{-1},\ldots,PA_mP^{-1})\in 
%           \utom(K)$ for some $P\in\pgltok$}\}. \]
In \cite{Nakamoto3}
we showed 
that there exists a scheme ${\rm Rep}_2(m)_{u}$
of finite type
over $\mathbb{Z}[1/2]$ (or $\mathbb{Z}/2{\mathbb{Z}}$).
Hence 
${\rm Rep}_2(m)_u(K)$ is an
algebraic variety over $K$
associated to the scheme ${\rm Rep}_2(m)_u$.
Note that 
there is a map
$\utom(K)\to {\rm Rep}_2(m)_{u}(K)$
of algebraic varieties
which is injective as a map 
of sets.

Let $\bto(K)$ be the subgroup 
of ${\rm PGL}_2(K)$ consisting
of upper triangular matrices.
The group
$\bto(K)$ acts on $\utom(K)$ by conjugation.
Since ${\rm PGL}_2(K)$ 
acts on ${\rm Rep}_2(m)_{u}(K)$ by conjugation,
the map $\utom(K)\to {\rm Rep}_2(m)_{u}(K)$
extends to a map ${\rm PGL}_2(K)\times \utom(K)\to {\rm Rep}_2(m)_{u}(K)$
of algebraic varieties.
This map factors through the quotient algebraic variety
${\rm PGL}_2(K)\times_{\bto(K)} \utom(K)$,
and hence
we obtain a map
\[ 
{\rm PGL}_2(K)\times_{\bto(K)}\utom(K)\longrightarrow {\rm Rep}_2(m)_u(K) 
\]
of algebraic varieties.

It is easy to prove the following lemma .

\begin{lemma}\label{lemma:normalizer-un}
If $P\in\pgltok$ satisfies $P\nto P^{-1}=\nto$,
then $P\in \bto(K)$.
\end{lemma}

\begin{theorem}
\label{theorem:description-repun}
The map
${\rm PGL}_2(K)\times_{\bto(K)}\utom(K)\to {\rm Rep}_2(m)_u(K)$
of algebraic varieties is a bijection 
\[ {\rm PGL}_2(K)\times_{\bto(K)}\utom(K)
   \stackrel{\cong}{\longrightarrow} \repuntok\]
as a map of sets.
\end{theorem}

\proof
For any $A=(A_1,\ldots,A_m)\in\repuntok$,
there exists $P\in\pgltok$ such that
$PAP^{-1}\in \nto$.
So 
$\pgltok\times_{\bto(K)} \utom(K)\to \repuntok$
is surjective.
Suppose that $PUP^{-1}=QVQ^{-1}$ for
$P,Q\in \pgltok$
and $U,V\in \utom(K)$.
Then $U=P^{-1}QVQ^{-1}P\in \utom(K)$.
So both of $U$ and $P^{-1}QVQ^{-1}P$
generate $\nto$. 
By Lemma~\ref{lemma:normalizer-un},
$P^{-1}Q=B\in \bto(K)$.
We see that
$(Q,V)=(PB,V)=(P,BVB^{-1})=(P,U)$
in $\pgltok\times_{\bto(K)}\utom(K)$.
This completes the proof.
\qed

\begin{remark}
\label{remark:description-repun}
\rm
\if0
There exists a scheme ${\rm Rep}_2(m)_{u}$
of finite type
over $\mathbb{Z}[1/2]$ (or $\mathbb{Z}/2{\mathbb{Z}}$)
by \cite{Nakamoto3},
and $\repuntok$ is the associated algebraic variety.
\fi
If the characteristic of $K$ is not $2$,
then the map
%\[ {\rm PGL}_2(K)\times_{\bto(K)}\utom(K)\stackrel{\cong}{\longrightarrow}
%   {\rm Rep}_2(m)_u(K) \]
is an isomorphism 
of algebraic varieties.
If the characteristic of $K$ is $2$,
then the map
induces a purely inseparable extension of degree $2$
between the function fields.
\end{remark}

We define the character variety ${\rm Ch}_2(m)_u(K)$ 
of unipotent representations of degree $2$
as the quotient algebraic variety
of ${\rm Rep}_2(m)_u(K)$
by ${\rm PGL}_2(K)$:
\[ 
{\rm Ch}_2(m)_u(K) = {\rm Rep}_2(m)_u(K)/{\rm PGL}_2(K).
\]
The map
${\rm PGL}_2(K)\times_{\bto(K)}\utom(K)\to {\rm Rep}_2(m)_u(K)$
induces a map
\[ \utom(K)/\bto(K)\longrightarrow {\rm Ch}_2(m)_u(K)\]
of algebraic varieties. 

Recall that for $(A_1,\ldots,A_m)\in ({\rm M}_2(K))^m$,
we have
$a_{ij}=(a_{ij}(1),\ldots,a_{ij}(m))$
for $i,j=1,2$,
where $A_k=(a_{ij}(k))\ (k=1,\ldots,m)$.
The map
$({\rm M}_2(K))^m\to (K^m)^2$
given by $(A_1,\ldots,A_m)\mapsto (a_{11},a_{12})$
induces an isomorphism
\[ \utom(K)\cong K^m\times (K^{m}-0).\]
of algebraic varieties.
By this isomorphism,
we obtain an isomorphism
\[ \utom(K)/\bto(K)\cong K^m\times {\rm P}^{m-1}(K),\]
of algebraic varieties.
This induces a map
\[ K^m\times {\rm P}^{m-1}(K)\longrightarrow
   {\rm Ch}_2(m)_u(K) \]
of algebraic varieties.

We easily obtain the following corollary 
by Theorem~\ref{theorem:description-repun}.

\begin{corollary}
\label{cor:descriptopn-chun}
The map
$K^m\times {\rm P}^{m-1}(K)\to {\rm Ch}_2(m)_u(K)$
of algebraic varieties is a bijection
\[ K^m\times {\rm P}^{m-1}(K)\stackrel{\cong}{\longrightarrow}
   \chuntok\]
as a map of sets.
\end{corollary}

\begin{remark}\rm
\label{remark:descriptopn-chun}
\if
There exists a scheme ${\rm Ch}_2(m)_{u}$
of finite type over $\mathbb{Z}[1/2]$ 
(or $\mathbb{Z}/2\mathbb{Z})$ by \cite{Nakamoto3},
and $\chuntok$ is the associated algebraic variety.
\fi
If the characteristic of $K$ is not $2$,
then the map 
%\[ K^m\times {\rm P}^{m-1}(K)
%   \stackrel{\cong}{\longrightarrow}{\rm Ch}_2(m)_u \]
is an isomorphism of algebraic varieties.
If the characteristic of $K$ is $2$,
then the map
induces a purely inseparable extension of degree $2$
between the function fields.
\end{remark}

\subsection{Cohomology groups of the moduli spaces of 
unipotent representations of degree $2$}

In this subsection we study the integral cohomology groups of
the moduli spaces of unipotent representations
of degree $2$ over $\mathbb{C}$.
First, we treat ${\rm Rep}_2(1)_{u}(\mathbb{C})$.

\begin{proposition}\label{prop:rep21un}
The space ${\rm Rep}_2(1)_{u}(\mathbb{C})$ is homotopy equivalent
to ${\rm P}^3(\mathbb{R})$.
Hence we have an isomorphism of commutative graded  algebras
\[ H^*({\rm Rep}_2(1)_{u}(\mathbb{C}))\cong \Lambda(s,t)/(2s,s^2,st)\]
with $|s|=2$ and $|t|=3$. 
\end{proposition}

\proof
By Theorem~\ref{theorem:description-repun}
and Remark~\ref{remark:description-repun},
${\rm Rep}_2(1)_{u}(\mathbb{C})\cong 
\pglto\times_{\bto(\mathbb{C})}{\rm U}_2(1)(\mathbb{C})$.
Recall that ${\rm U}_2(1)\cong \mathbb{C}\times \mathbb{C}^{\times}$.
Then $\bto(\mathbb{C})$ 
acts trivially on the left factor $\mathbb{C}$
and transitively on the right factor $\mathbb{C}^{\times}$.
Let $S$ be the stabilizer subgroup of $\bto(\mathbb{C})$
at $1\in \mathbb{C}^{\times}$.
Then 
\[  S=\left\{\left(\begin{array}{cc}
                      1 & * \\
                      0 & 1 \\
                   \end{array}\right)\in\pglto\right\}.\]
We can write ${\rm U}_2(1)(\mathbb{C})
\cong \mathbb{C}\times (\mathbb{}\bto(\mathbb{C})/S)$.
Hence ${\rm Rep}_2(1)_{u}(\mathbb{C})
\cong (\pglto/S)\times \mathbb{C}\simeq \pglto/S$.
Since $S\cong \mathbb{C}$,
$\pglto\to \pglto/S$ is a homotopy equivalence.
So we see that ${\rm Rep}_2(1)_u(\mathbb{C})\simeq \pglto$.
It is well-known that the inclusion ${\rm PU}(2)\to \pglto$
induces a homotopy equivalence
${\rm PU}(2)\simeq \pglto$.
Hence 
${\rm Rep}_2(1)_{u}(\mathbb{C})$ is homotopy equivalent
to ${\rm PU(2)}\cong {\rm P}^3(\mathbb{R})$.
\qed

\begin{corollary}
We have an isomorphism of commutative graded algebras
\[ H^*({\rm Rep}_2(1)_{u}(\mathbb{C});\mathbb{Q})\cong \Lambda(t)\]
with $|t|=3$.
\end{corollary}

\bigskip

For $m\ge 2$ we have the following theorem on the
cohomology groups of $\repunto$.

%By Porposition~\ref{prop:rep21un},
%we can calculate the cohomology of
%${\rm Rep}_2(1)_{un}$.

\begin{theorem}\label{theorem:integral-cohomology-repunto-mge2}
For $m\ge 2$, we have an isomorphism of commutative graded  algebras
\[ H^*(\repunto)\cong \Lambda (u,s)/(u^2) \]
with $|u|=2$ and $|s|=2m-1$.
\end{theorem}

\proof
By Theorem~\ref{theorem:description-repun}
and Remark~\ref{remark:description-repun},
there is a fibre bundle
$\utom(\mathbb{C})\to
\repunto\to\pglto/\bto(\mathbb{C})
\cong {\rm P}^1(\mathbb{C})$.
Since $\utom(\mathbb{C})\simeq S^{2m-1}$,
the associated Serre spectral sequence collapses.
This completes the proof.
\qed

\begin{corollary}
For $m\ge 2$, we have an isomorphism of 
commutative graded algebras
\[ H^*(\repunto;\mathbb{Q})\cong \Lambda (u,s)/(u^2) \]
with $|u|=2$ and $|s|=2m-1$.
\end{corollary}

\proof
This follows from the fact that the integral cohomology groups
of $\repunto$ is torsion-free 
by Theorem~\ref{theorem:integral-cohomology-repunto-mge2}.
\qed

\

By Corollary~\ref{cor:descriptopn-chun}
and Remark~\ref{remark:descriptopn-chun},
the cohomology groups of $\chunto$
is given as follows.

\begin{proposition}
We have isomorphisms of commutative graded algebras
\[ H^*(\chunto)\cong \mathbb{Z}[t]/(t^m)\]
and
\[ H^*(\chunto;\mathbb{Q})\cong \mathbb{Q}[t]/(t^m)\]
with $|t|=2$.
\end{proposition}

%\bigskip
%
%We can also calculate the cohomology
%of $\chunto$ by Corollary~\ref{cor:descriptopn-chun}.

\if0
\begin{remark}\rm
By Corollary~\ref{cor:descriptopn-chun},
Proposition~\ref{prop:rep21un},
and Theorem~\ref{theorem:integral-cohomology-repunto-mge2},
we can easily calculate 
the rational cohomology groups of
$\repunto$ and $\chunto$.
\end{remark}
\fi

%\input{sc}
\section{The moduli spaces of scalar representations}
\label{section:scalar-representations}

In this section we define the moduli spaces related to 
scalar representations.
It is easy to give descriptions of the moduli spaces.
Then we obtain the cohomology groups of them.

Let $K$ be an algebraically closed field.
We define the moduli space ${\rm Rep}_n(m)_{sc}(K)$ 
of scalar representations by
\[ {\rm Rep}_n(m)_{sc}(K)=\{(A_1,\ldots,A_m)\in (M_n(K))^m|\ 
    \dim_K K\langle A_1,\ldots,A_m\rangle=1\},\]
where $K\langle A_1,\ldots,A_m\rangle$
is the $K$-subalgebra of $M_n(K)$ generated by
$A_1,\ldots,A_m$.

\begin{theorem}\label{thm:descriptionRepsc}
We have an isomorphism of smooth algebraic varieties 
\[ {\rm Rep}_n(m)_{sc}(K)\cong K^m. \]
Hence we have isomorphisms of commutative graded algebras
\[ H^*({\rm Rep}_n(m)_{sc}(\mathbb{C}))\cong \mathbb{Z}\]
and
\[ H^*({\rm Rep}_n(m)_{sc}(\mathbb{C});\mathbb{Q})\cong \mathbb{Q}.\]
\end{theorem}

\proof
If $(A_1,\ldots,A_m)\in {\rm Rep}_n(m)_{sc}(K)$,
then $A_i$ is a scalar matrix for $1\le i\le m$.
Hence ${\rm Rep}_n(m)_{sc}(K)\cong K^m$.
\qed

\bigskip

The group $\pglk$ acts on
${\rm Rep}_n(m)_{sc}(K)$ by conjugation.
We define the character variety ${\rm Ch}_n(m)_{sc}(K)$
of scalar representations by
the quotient space
\[ {\rm Ch}_n(m)_{sc}(K)= {\rm Rep}_n(m)_{sc}(K)/\pglk.\]

\begin{theorem}\label{thm:descriptionChsc}
We have an isomorphism of smooth algebraic varieties
\[ {\rm Ch}_n(m)_{sc}(K)\cong K^m. \]
Hence we have isomorphisms of commutative graded algebras
\[ H^*({\rm Ch}_n(m)_{sc}(\mathbb{C}))\cong \mathbb{Z}\]
and
\[ H^*({\rm Ch}_n(m)_{sc}(\mathbb{C});\mathbb{Q})\cong \mathbb{Q}\]
\end{theorem}

\proof
This follows from the fact that the action of $\pglk$ on
${\rm Rep}_n(m)_{sc}(K)$ is trivial 
since ${\rm Rep}_n(m)_{sc}(K)$ consists of
$m$-tuples of scalar matrices.
\qed

\begin{remark}\rm
There exist smooth schemes 
${\rm Rep}_n(m)_{sc}\cong \mathbb{A}^m_{\mathbb{Z}}$
and 
${\rm Ch}_n(m)_{sc}\cong \mathbb{A}^m_{\mathbb{Z}}$
over $\mathbb{Z}$,
and ${\rm Rep}_n(m)_{sc}(K)$ and
${\rm Ch}_n(m)_{sc}(K)$ are the 
associated algebraic varieties.
\end{remark}

\section{Virtual Hodge polynomials of the moduli spaces}
\label{section:VHP-air}

In this section we study the virtual Hodge polynomials of 
the moduli spaces of representations of degree $2$
over $\mathbb{C}$.
See, for example, 
%\cite{Cheah1, Cheah2, Danilov}
\cite{Danilov}
for the precise definition and properties
of the virtual Hodge polynomial.
Also, see \cite{topos,rational} for
the virtual Hodge polynomials of the moduli spaces
of representations with Borel mold.

For a mixed Hodge structure $(V,W,F)$,
we denote by $a^{p,q}(V)$
the dimension of the $(p,q)$-component of the pure Hodge
structure ${\rm Gr}_{p+q}^W(V)$ of weight $p+q$.
For an algebraic scheme $X$ over $\mathbb{C}$,
we denote by ${\rm VHP}(X)$ the virtual Hodge polynomial of $X$:
\[ {\rm VHP}(X):=\sum_{p,q,n}(-1)^na^{p,q}(H^n(X;\mathbb{Q}))x^py^q.\]
We also denote by ${\rm VHP}_c(X)$ the virtual Hodge polynomial
of $X$ based on the compact support cohomology.
Note that if $X$ is smooth of pure dimension $m$,
then 
\begin{align}\label{equation:poincare-duality} 
{\rm VHP}_c(X)=(xy)^m{\rm VHP}(X)
   \left(x^{-1},y^{-1}\right) 
\end{align}
by the Poincar\'e duality.
For simplicity, we set $z=xy$.

\subsection{Virtual Hodge polynomials of the moduli spaces of
  representations of degree $2$}

In this subsection we study the virtual Hodge polynomials of
the moduli spaces of representations of degree $2$
over $\mathbb{C}$.
In particular, we calculate the virtual
Hodge polynomial of the moduli space of absolutely irreducible
representations of degree $2$.

Let $K$ be an algebraically closed field.
Let ${\rm Rep}_n(m)(K)=({\rm M}_n(K))^m$
be the representation variety of degree $n$
for the free monoid with $m$ generators. 
We define the subvariety ${\rm Rep}_n(m)_{{\rm rk}h}(K)$ of
${\rm Rep}_n(m)(K)$ by
\[ {\rm Rep}_n(m)_{{\rm rk}h}(K)=\{(A_1,\ldots,A_m)\in ({\rm
  M}_n(K)^m |\ 
    \dim_K K\langle A_1,\ldots,A_m\rangle=h\}.\]
The representation variety $\repairk$
of absolutely irreducible representations
is defined by
\[ \repairk = {\rm Rep}_n(m)_{{\rm rk}n^2}(K).\]

\begin{remark}\rm
By \cite{Nakamoto1},
there exists a smooth scheme 
${\rm Rep}_n(m)_{air}$ over $\mathbb{Z}$ 
and 
$\repairk$ is the associated algebraic variety.
\end{remark}

Let 
${\rm Rep}_n(m)_B(K)$ be the representation variety 
of representations with Borel mold.
When $n=2$,
we have ${\rm Rep}_2(m)_B(K)={\rm Rep}_2(m)_{{\rm rk}3}(K)$.
We calculated 
the virtual Hodge polynomial 
of ${\rm Rep}_n(m)_B(\mathbb{C})$ in 
\cite{topos} and \cite{rational}.
%When $n=2$, the virtual Hodge polynomial of
%${\rm Rep}_2(m)_B$ is given as follows.

%\newpage

\begin{proposition}[{\cite[Proposition~7.9]{topos}} and
{\cite[Corollary~8.16]{rational}}]
\label{prop:VHP_B}
The virtual Hodge polynomial of ${\rm Rep}_n(m)_B(\mathbb{C})$ 
is given by
%\[ \begin{array}{rcl}
%%     {\rm VHP}({\rm Rep}_2(m)_B)&=&(1-z^m)(1-z^{m-1})(1+z),\\[3mm]
%%     {\rm VHP}_c({\rm Rep}_2(m)_B)&=& z^{m+1}(z^m-1)(z^{m-1}-1)(z+1).\\
%     {\rm VHP}({\rm Rep}_n(m)_B(\mathbb{C}))&=&
%     \prod_{k=1}^{n-1}(1-kz^m)\cdot 
%     \left( \frac{1-z^{m-1}}{1-z}\right)^{n-1}\cdot
%     \prod_{i=2}^n(1-z^i),\\[3mm]
%     {\rm VHP}_c({\rm Rep}_n(m)_B(\mathbb{C}))&=& 
%     \frac{(z^m-z)^{n-1}}{(z-1)^{n-1}}\cdot z^{m(n-2)(n-1)/2}\cdot
%     \prod_{k=0}^{n-1}(z^m-k)\cdot \prod_{k=2}^n (z^k-1).\\
%  \end{array}\]
\[ \begin{array}{rcl}
     {\rm VHP}({\rm Rep}_n(m)_B(\mathbb{C}))&=&\displaystyle
     \frac{(1-z^{m-1})^{n-1}\prod_{k=1}^{n-1}(1-kz^m)
     \prod_{i=1}^n(1-z^i)}
     {(1-z)^n},\\[3mm]
     {\rm VHP}_c({\rm Rep}_n(m)_B(\mathbb{C}))&=&\displaystyle 
     \frac{z^{m(n-1)(n-2)/2}(z^m-z)^{n-1}
     \prod_{k=0}^{n-1}(z^m-k)\prod_{k=1}^n (z^k-1)}
     {(z-1)^n}.\\
  \end{array}\]
\end{proposition}

\

Next we consider the virtual Hodge polynomial
of $\reptoss$.

\begin{proposition}\label{prop:VHP_reptoss}
The virtual Hodge polynomial of $\reptoss$ is given by
\[ \begin{array}{rcl}
     {\rm VHP}(\reptoss)&=&1-z^m,\\[3mm]
     {\rm VHP}_c(\reptoss)&=&z^{m+2}(z^m-1).\\
   \end{array}\]
\end{proposition}

\proof
By Corollary~\ref{cor:rationaliso-config-repss},
the inclusion $\configto\hookrightarrow \reptoss$
induces an isomorphism of rational cohomology groups.
Hence 
${\rm VHP}(\reptoss)={\rm VHP}(\configto)=1-z^m$.
Since $\dim_{\mathbb{C}}\reptoss=2m+2$,
we have
\[ {\rm VHP}_c(\reptoss)(z)=z^{2m+2}{\rm VHP}(\reptoss)(z^{-1}). \]
Hence we obtain that
${\rm VHP}_c(\reptoss)=z^{m+2}(z^m-1)$.
\qed

\

The virtual Hodge polynomial of $\repunto$
is given as follows.

\begin{proposition}\label{prop:VHP_repunto}
The virtual Hodge polynomial of $\repunto$ is given by
\[ \begin{array}{rcl}
     {\rm VHP}(\repunto)&=&(1+z)(1-z^m),\\[3mm]
     {\rm VHP}_c(\repunto)&=&z^m(z+1)(z^m-1).\\
   \end{array}\]
\end{proposition}

\proof
By Theorem~\ref{theorem:description-repun}
and Remark~\ref{remark:description-repun},
we have the fiber bundle  
$\utom(\mathbb{C})\to \repunto\to 
\pglto/\bto(\mathbb{C})$ 
with respect to the Zariski topology.
Note that $\utom(\mathbb{C})
\cong \mathbb{C}^m\times (\mathbb{C}^m-0)$
and $\pglto/\bto(\mathbb{C})
\cong {\rm P}^1(\mathbb{C})$.
By the property of the virtual Hodge polynomial,
we obtain that 
\[ {\rm VHP}_c(\repunto)={\rm VHP}_c(\utom(\mathbb{C}))\cdot 
{\rm VHP}_c(\pglto/\bto(\mathbb{C})). \]
Since ${\rm VHP}_c(\utom(\mathbb{C}))=z^m(z^m-1)$
and ${\rm VHP}_c(\pglto/\bto(\mathbb{C}))=1+z$,
%we see that
${\rm VHP}_c(\repunto)=z^m(z+1)(z^m-1)$.
Since $\dim_{\mathbb{C}}\repunto=2m+1$,
we have
\[ {\rm VHP}(\repunto)=z^{2m+1}{\rm VHP}_c(\repunto)(z^{-1}). \]
Hence ${\rm VHP}(\repunto)=(1+z)(1-z^m)$.
%Similarly, we obtain that
%${\rm VHP}(\repunto)={\rm VHP}(\utom)\cdot
%{\rm VHP}(\pgl/\bto)=(1-z^m)(1+z)$.
\qed

\bigskip

%We define ${\rm Rep}_2(m)_{{\rm rk}2}$ by
%\[ {\rm Rep}_2(m)_{{\rm rk}2}=\{(A_1,\ldots,A_m)\in ({\rm
%  M}_2(\mathbb{C}))^m |\ 
%    \dim\mathbb{C}\langle A_1,\ldots,A_m\rangle=2\}.\]
We consider the virtual Hodge polynomial
of ${\rm Rep}_2(m)_{{\rm rk}2}(\mathbb{C})$.
Let $U_i$ be the subspace of ${\rm Rep}_2(m)_{{\rm rk}2}(\mathbb{C})$
consisting of $(A_1,\ldots,A_m)$ such that
$A_1,\ldots,A_{i-1}$ are scalar matrices, and
$A_i$ is not a scalar matrix.
Then we have a decomposition of ${\rm Rep}_2(m)_{{\rm rk}2}(\mathbb{C})$:
\[ {\rm Rep}_2(m)_{{\rm rk}2}(\mathbb{C})=U_1\cup\cdots \cup U_m.\]
Note that $U_i\cap U_j=\emptyset$ for $i\neq j$.
Since
$U_k$ is open in $\cup_{i=k}^nU_i$, 
we have 
\[ {\rm VHP}_c(\cup_{i=k}^nU_i)={\rm VHP}_c(U_k)+
{\rm VHP}_c(\cup_{i=k+1}^nU_i) \]
for $1\le k\le n$. 
Hence we obtain 
\[ {\rm VHP}_c({\rm Rep}_2(m)_{{\rm rk}2}(\mathbb{C}))
   =\sum_{i=1}^n{\rm VHP}_c(U_i).\]
Let $I_2$ be the identity matrix in ${\rm M}_2(\mathbb{C})$.
For $(A_1,\ldots,A_m)\in U_i$,
we can uniquely write $A_r=\alpha_rI_2+\beta_rA_i$
for $r\ge i+1$, where $\alpha_r,\beta_r\in\mathbb{C}$.
This implies that
\[ U_i\cong \mathbb{C}^{i-1}\times ({\rm M}_2(\mathbb{C})-\mathbb{C}
            \cdot I_2)\times(\mathbb{C}^2)^{m-i}.\]    
Then we have
${\rm VHP}(U_i)=1-z^3$ and ${\rm VHP}_c(U_i)=z^{i-1}\cdot (z^4-z)\cdot
z^{2(m-i)}=z^{2m-i}(z^3-1)$.
Hence we obtain the following proposition.

\begin{proposition}\label{prop:VHPcRepRk2}
The virtual Hodge polynomial of ${\rm Rep}_2(m)_{{\rm rk}2}(\mathbb{C})$ is 
given by
\[ {\rm VHP}_c({\rm Rep}_2(m)_{{\rm rk}2}(\mathbb{C}))=
   z^m(z^2+z+1)(z^m-1).\]
\end{proposition}

\if0
\proof
This follows from the fact that 
${\rm VHP}_c({\rm Rep}_2(m)_{{\rm rk}2})=\sum_{i=1}^m{\rm VHP}_c(U_i)$.
\qed
\fi

\begin{remark}\rm
We have a decomposition 
\[ {\rm Rep}_2(m)_{{\rm rk}2}(\mathbb{C})=\reptoss\cup \repunto, \]
where $\reptoss\cap \repunto=\emptyset$
and
$\reptoss$ is open in ${\rm Rep}_2(m)_{{\rm rk}2}(\mathbb{C})$.
By Propositions~\ref{prop:VHP_reptoss}
and \ref{prop:VHP_repunto},
${\rm VHP}_c(\reptoss)=z^{m+2}(z^m-1)$ and
${\rm VHP}_c(\repunto)=z^m(z+1)(z^m-1)$.
Hence we have
\[ \begin{array}{rcl}
     {\rm VHP}_c({\rm Rep}_2(m)_{{\rm rk}2}(\mathbb{C}))&=&
      z^{m+2}(z^m-1)+z^m(z+1)(z^m-1)\\[3mm]
     &=&z^m(z^2+z+1)(z^m-1). \\
   \end{array}\]
\end{remark}

By definition,
we have ${\rm Rep}_2(m)_{{\rm rk}1}(\mathbb{C})=
{\rm Rep}_2(m)_{{\rm sc}}(\mathbb{C})$.
Then 
we obtain the following proposition
by Theorem~\ref{thm:descriptionRepsc}.

\begin{proposition}\label{prop:VHP_sc}
The virtual Hodge polynomial of ${\rm Rep}_n(m)_{sc}(\mathbb{C})$ 
is given by
\[ \begin{array}{rcl}
 {\rm VHP}({\rm Rep}_n(m)_{sc}(\mathbb{C}))&=&1,\\[3mm]
 {\rm VHP}_c({\rm Rep}_n(m)_{sc}(\mathbb{C}))&=&z^m.\\
   \end{array}\]
\end{proposition}

%We say that a representation
%$\Upsilon_m\to M_n(\mathbb{C})$ for
%the free monoid $\Upsilon_m$ of rank $m$ 
%is absolutely irreducible if
%the image of the generators generates
%the full matrix ring ${\rm M}_n(\mathbb{C})$
%as a $\mathbb{C}$-algebra.
%
%We define $\repair$ by
%\[ \repair=\{(A_1,\ldots,A_m)\in ({\rm M}_n(\mathbb{C}))^m|\ 
%   \dim\mathbb{C}
%   \langle A_1,\ldots,A_m\rangle =n^2\}.\]
%Then $\repair$ is the moduli space of absolutely irreducible
%representations for $\Upsilon_m$ of degree $n$. 

Recall that $\repairk={\rm Rep}_n(m)_{{\rm rk}n^2}(K)$
is the representation variety of absolutely irreducible representations. 
Since $\repairk$ is an open subvariety of 
${\rm Rep}_n(m)(K)=({\rm M}_n(K))^m$,
$\repairk$ is smooth of pure dimension $mn^2$.
When $n=2$, we can calculate the virtual Hodge polynomial 
of $\reptoair$ by using the above results.

\begin{theorem}[Theorem~\ref{theorem:Main-Theorem1}]\label{thm:VHPreptoair}
The virtual Hodge polynomial of $\reptoair$
is given by
\[  \begin{array}{rcl}
     {\rm VHP}(\reptoair)&=&(1-z^m)(1-z^{m-1}),\\[3mm]
     {\rm VHP}_c(\reptoair)&=&z^{2m+1}(z^m-1)(z^{m-1}-1).\\
    \end{array}\]
\end{theorem}

\proof
First, we calculate ${\rm VHP}_c(\reptoair)$.
We see that
${\rm VHP}_c({\rm Rep}_2(m))=z^{4m}$
since ${\rm Rep}_2(m)(\mathbb{C})={\rm M}_2(\mathbb{C})^m$.
We have a decomposition of ${\rm Rep}_2(m)(\mathbb{C})$: 
\[ {\rm Rep}_2(m)(\mathbb{C})=\cup_{h=1}^4 
   {\rm Rep}_2(m)_{{\rm rk}h}(\mathbb{C}).\]
%where 
%\[ {\rm Rep}_2(m)_{{\rm rk}h}=
%   \{(A_1,\ldots,A_m)\in ({\rm M}_2(\mathbb{C}))^m|\ 
%   \dim\mathbb{C}
%   \langle A_1,\ldots,A_m\rangle =h\}.\]
%\[ {\rm Rep}_2(m)=\reptoair\cup {\rm Rep}_2(m)_B\cup
%   \repunto\cup \reptoss\cup {\rm Rep}_2(m)_{sc}.\]
Then ${\rm Rep}_2(m)_{{\rm rk}i}(\mathbb{C})\cap 
      {\rm Rep}_2(m)_{{\rm rk}j}(\mathbb{C})=\emptyset$ if $i\neq j$.
Furthermore,
${\rm Rep}_2(m)_{{\rm rk}h}(\mathbb{C})$ is closed in 
$\cup_{i=h}^4 {\rm Rep}_2(m)_{{\rm rk}i}(\mathbb{C})$.
By the additivity property of virtual Hodge polynomial,
\[ {\rm VHP}_c({\rm Rep}_2(m)(\mathbb{C}))
   =
   \sum_{h=1}^4\ {\rm VHP}_c({\rm Rep}_2(m)_{{\rm rk}h}(\mathbb{C})).\]
Recall that ${\rm Rep}_2(m)_{{\rm rk}4}(\mathbb{C})=\reptoair$ and
${\rm Rep}_2(m)_{{\rm rk}3}(\mathbb{C})={\rm Rep}_2(m)_B(\mathbb{C})$.
Then we can calculate ${\rm VHP}_c(\reptoair)$ by
Propositions~\ref{prop:VHP_reptoss},
\ref{prop:VHP_repunto},
\ref{prop:VHP_sc}, and
\ref{prop:VHP_B}.

Since $\reptoair$ is smooth of pure dimension $4m$, 
we can calculate ${\rm VHP}(\reptoair)$
by the formula ${\rm VHP}(\reptoair)
=z^{4m}{\rm VHP}_c(\reptoair)(z^{-1})$.
\qed

\subsection{Virtual Hodge polynomials of the character varieties
of degree $2$}

In this subsection we study the virtual Hodge polynomials of
the character varieties of degree $2$ over $\mathbb{C}$.

%\subsection{Virtual Hodge polynomials of the character varieties
%of degree $2$}
%
%In this section
%we calculate the  virtual Hodge polynomial of
%the character varieties of representations of degree $2$. 

\begin{lemma}
The virtual Hodge polynomial of the configuration space
$\unconfig$ is given by
\[ \begin{array}{rcl}
    {\rm VHP}(\unconfig)&=&1-z^m,\\[3mm]
    {\rm VHP}_c(\unconfig)&=&z^{m(n-1)}(z^m-1).\\
   \end{array}\]    
\end{lemma}

\proof
By Lemma~\ref{lemma:cohomology-unconfig},
the quotient map $\config\to \unconfig$ induces
an injection on the rational cohomology groups
$H^*(\unconfig;\mathbb{Q})\hookrightarrow 
H^*(\config;\mathbb{Q})$,
and $H^*(\unconfig;\mathbb{Q})=\Lambda(s)$ with $|s|=2m-1$.
Since the mixed Hodge structure on $H^{2m-1}(\config;\mathbb{Q})$
is pure of type $(m,m)$ (cf.~\cite[Section~5.3]{rational}),
we see that the virtual Hodge polynomial of $\unconfig$
is given by ${\rm VHP}(\unconfig)=1-z^m$.
Since $\unconfig$ is smooth of dimension $mn$,
${\rm VHP}_c(\unconfig)=z^{mn}{\rm VHP}(\unconfig)(z^{-1})$.
Hence we obtain that 
${\rm VHP}_c(\unconfig)=z^{m(n-1)}(z^m-1)$.
\qed

\

By Corollary~\ref{cor:chss-unconfig},
we obtain the following proposition.

\begin{proposition}
The virtual Hodge polynomial of $\chss$ is given by
\[ \begin{array}{rcl}
    {\rm VHP}(\chss)&=&1-z^m,\\[3mm]
    {\rm VHP}_c(\chss)&=&z^{m(n-1)}(z^m-1).\\
   \end{array}\]    
\if0
In particular, we have
\[ \begin{array}{rcl}
    {\rm VHP}({\rm Ch}_2(m)_{ss})&=&1-z^m,\\[3mm]
    {\rm VHP}_c({\rm Ch}_2(m)_{ss})&=&z^{m}(z^m-1).\\
   \end{array}\]    
\fi
\end{proposition}

\if0
\proof
The composition $\config\to \repss\to \chss$
is a map of algebraic schemes over $\mathbb{C}$.
By Proposition~\ref{prop:chss-rational},
the mixed Hodge structure on
$H^{2m-1}(\chss)$ is a subobject of $H^{2m-1}(\config)$.
The mixed Hodge structure on $H^{2m-1}(\config)$ is
homogeneous, 
we obtain that
${\rm VHP}(\chss)=1-z^m$.
Since $\dim \chss=mn$,
we obtain that
${\rm VHP}_c(\chss)=z^{mn}(1-1/z^m)
=z^{m(n-1)}(z^m-1)$.
\qed
\fi

By Corollary~\ref{cor:descriptopn-chun}
and Remark~\ref{remark:descriptopn-chun},
we obtain the following proposition. 

\begin{proposition}
The virtual Hodge polynomial of $\chunto$ is given by
\[ \begin{array}{rcl}
     {\rm VHP}(\chunto)&=&
     {\displaystyle \frac{1-z^m}{1-z}},\\[3mm]
     {\rm VHP}_c(\chunto)&=&
     {\displaystyle  \frac{z^m(z^m-1)}{z-1}}.\\
   \end{array}\]
\end{proposition}

\bigskip

By Theorems~\ref{thm:descriptionRepsc}
and \ref{thm:descriptionChsc},
we have ${\rm Ch}_n(m)_{sc}(\mathbb{C})=
{\rm Rep}_n(m)_{sc}(\mathbb{C})\cong 
\mathbb{C}^m$.
Hence we obtain the following proposition.
%by Proposition~\ref{prop:VHP_sc}.

\begin{proposition}
The virtual Hodge polynomial of ${\rm Ch}_n(m)_{sc}(\mathbb{C})$
is given by
\[ \begin{array}{rcl}
 {\rm VHP}({\rm Ch}_n(m)_{sc}(\mathbb{C}))&=&1,\\[3mm]
 {\rm VHP}_c({\rm Ch}_n(m)_{sc}(\mathbb{C}))&=&z^m.\\
   \end{array}\]
\if0
In particular, we have
\[ \begin{array}{rcl}
 {\rm VHP}({\rm Ch}_2(m)_{sc})&=&1,\\[3mm]
 {\rm VHP}_c({\rm Ch}_2(m)_{sc})&=&z^m.\\
   \end{array}\]
\fi
\end{proposition}

We calculated 
the virtual Hodge polynomial 
of ${\rm Ch}_n(m)_B(\mathbb{C})$ in 
\cite{topos} and \cite{rational}.

\begin{proposition}[{\cite[Proposition~7.8]{topos} and 
\cite[Corollary~8.8]{rational}}]
The virtual Hodge polynomial of
${\rm Ch}_n(m)_B(\mathbb{C})$ is given by
%\[ \begin{array}{rcl}
%     {\rm VHP}({\rm Ch}_n(m)_B(\mathbb{C}))&=&
%     {\displaystyle \prod_{k=1}^{n-1}(1-k z^m)\cdot 
%       \left(\frac{1-z^{m-1}}{1-z}\right)^{n-1} },\\[3mm]
%     {\rm VHP}_c({\rm Ch}_n(m)_B(\mathbb{C})) &=& 
%     {\displaystyle \frac{(z^{m-1}-1)^{n-1}}{(z-1)^{n-1}}
%      z^{(m-1)(n-2)(n-1)/2}
%      \prod_{k=0}^{n-1} (z^m - k)}.\\ 
%   \end{array}\]
\[ \begin{array}{rcl}
     {\rm VHP}({\rm Ch}_n(m)_B(\mathbb{C}))&=&
     {\displaystyle \frac{(1-z^{m-1})^{n-1}
        \prod_{k=1}^{n-1}(1-k z^m)}
       {(1-z)^{n-1}}},\\[3mm]
     {\rm VHP}_c({\rm Ch}_n(m)_B(\mathbb{C})) &=& 
     {\displaystyle \frac{z^{(m-1)(n-1)(n-2)/2}
       (z^{m-1}-1)^{n-1}\prod_{k=0}^{n-1} (z^m - k)}
       {(z-1)^{n-1}}}.\\ 
   \end{array}\]
\end{proposition}

\bigskip

The conjugate action of matrices
induces an action of   
$\pgl$ on $\repair$.
The character variety $\chair$
of absolutely irreducible representations
is defined to be the quotient
space 
\[ \chair=\repair/\pgl. \]
We let $\pi$ be the quotient map
\[ \pi: \repair \to \chair. \]
%\begin{remark}\rm
By \cite{Nakamoto1},
there exists a smooth scheme 
${\rm Ch}_n(m)_{air}$ over $\mathbb{Z}$,
and $\chair$ is the associated algebraic variety.
Furthermore,
the quotient map $\pi$ is induced by
a map
\[ \pi: {\rm Rep}_n(m)_{air}\to {\rm Ch}_n(m)_{air} \]
of schemes over $\mathbb{Z}$.
%\end{remark}
%
%Let us verify if 
%$\pi$ satisfies the conidtions in
%Theorem~\ref{thm:product-decomposition-VHP-general}
%Theorem~\ref{thm:VHP-product-complex-topology}.
%We consider the Leray spectral sequence of $f$:
%\[ E_2^{p,q}=H^p(\chair;R^qf_*\underline{\mathbb{Q}})\Longrightarrow 
%           H^{p+q}(\repair).\]
%By \cite[??]{Nakamoto3},
%$f$ is a fibre bundle with respect to the \'etale topology.
%In particular,
%this implies that $f$ is a fibre bundle with repset to
%the classical topology.
%the complex topology.

To calculate the virtual Hodge polynomial of
$\chair$,
we need the multiplicative property
of the virtual Hodge polynomials.
Let $\varphi: X\to Y$ be a map 
of complex algebraic varieties. 
%We recall the relationship 
%between the Hodge polynomials of $X$ and $Y$
%under some conditions.
\if0
Then we have the Leray spectral sequence
\[ E_2^{p,q}=H^p(Y;R^q\varphi_*\underline{\mathbb{Q}}_X)\Longrightarrow 
           H^{p+q}(X),\]
where $R^q\varphi_*\underline{\mathbb{Q}}_X$ 
is the $q$th higher direct image
of the constant sheaf $\underline{\mathbb{Q}}_X$ over $X$.
Using the cohomology with compact support,
we have a compact support version of the Leray
spectral sequence
\[ E_2^{p,q}=H^p_c(Y;R^q\varphi_!\underline{\mathbb{Q}}_X)\Longrightarrow 
           H^{p+q}_c(X),\]
In this subsection we recall the relationship 
between the Hodge polynomials of $X$ and $Y$
under some conditions.
\fi
Let $\underline{\mathbb{Q}}_X$ be the constant sheaf on $X$.
We denote by $R^q\varphi_*\underline{\mathbb{Q}}_X$ 
the $q$th higher direct image
of $\underline{\mathbb{Q}}_X$,
and by $R^q\varphi_!\underline{\mathbb{Q}}_X$ 
the $q$th higher direct image with compact support.

\begin{theorem}
[{cf.~\cite[Lemma~2 and Remark~2]{CLMS} and 
\cite[Theorem~6.1]{Dimca-Lehrer}}]
%, and
%\cite[Theorem~6.5 and Corollary~14.14]{Peters-Steenbrink}}]
\label{thm:VHP-product-complex-topology}
Let $f: X\to Y$ be a map of complex algebraic varieties,
where $Y$ is smooth and connected.
We suppose that
$f$ is a locally trivial fibration with respect 
to the complex topology.
Suppose further that 
$R^q\varphi_*\underline{\mathbb{Q}}_X$,
respectively
$R^q\varphi_{!}\underline{\mathbb{Q}}_X$,  
is a constant sheaf on $Y$
for all $q$.
Then we have
\[ {\rm VHP}(X)  =  {\rm VHP}(Y)\cdot {\rm VHP}(F),\]
respectively
\[ {\rm VHP}_c(X)  =  {\rm VHP}_c(Y)\cdot {\rm VHP}_c(F).\]
\end{theorem}

\if0
\begin{lemma}
[{cf.~\cite[Lemma~2]{CLMS}, \cite[Theorem~6.1]{Dimca-Lehrer}, and
\cite[Theorem~6.5 and Corollary~14.14]{Peters-Steenbrink}}]
Let $f: X\to Y$ be a map of complex algebraic varieties,
where $Y$ is smooth and connected.
We suppose that
$f$ is a locally trivial fibration with respect 
to the complex topology.
Suppose further that 
$R^q\varphi_{!}\underline{\mathbb{Q}}_X$ 
is a constant sheaf on $Y$ for all $q$.
Then we have
\[ {\rm VHP}_c(X)  =  {\rm VHP}_c(Y)\cdot {\rm VHP}_c(F).\]
\end{lemma}
\fi

Let us verify if 
$\pi$ satisfies the conditions in
%Theorem~\ref{thm:product-decomposition-VHP-general}
Theorem~\ref{thm:VHP-product-complex-topology}.

\begin{proposition}\label{prop:loctrietale} 
Let $f : X \to Y$ be a principal fibre bundle with 
group $G$ over a scheme $S$
in the sense of \cite[Definition~0.10]{GIT}. 
In other words, $(Y, f)$ is a geometric quotient of $X$ by 
$G$ over $S$ satisfying 
\begin{enumerate}
\item $G$ is flat and of finite type over $S$,  
\item $f$ is a flat morphism of finite type,  
\item $G\times_{S} X \to X\times_{Y}X$ is an isomorphism. 
\end{enumerate} 
If $f$ is smooth, then $f$ has a local trivialization with 
respect to the \'{e}tale topology.  
\end{proposition}

\proof
By \cite[17.16.3(ii)]{EGAIV}
(cf. \cite[Chapter~I, Propositions~3.24 and 3.26]{Milne}), 
there exist  
a surjective \'{e}tale morphism $h: Y'\to Y$,
and a $Y$-morphism $g: Y'\to X$.
Hence $f$ has a local trivialization with respect to 
the \'{e}tale topology. 
\qed

\bigskip 

By \cite[Corollary~6.4]{Nakamoto1}, $\pi$ 
is a principal fibre bundle with ${\rm PGL}_n$. 
Note that $\pi$ is smooth
because $\pi$ is flat and the $\pi^{-1}(\bar{x})$
is regular for any geometric point $\bar{x}$
of ${\rm Ch}_n(m)_{air}$. 
By Proposition~\ref{prop:loctrietale},
%$\pi$ has local trivialization with 
%respect to the \'{e}tale topology.  
%Hence 
we obtain the following corollary. 

\begin{corollary}
The map $\pi$ is a fibre bundle with respect to 
the \'{e}tale topology.
In particular,
%$\pi$ is a fibre bundle 
so is $\pi$
with respect to the complex topology. 
\end{corollary}

Notice that $\chair$ is path-connected.

\begin{lemma}\label{lemma:triviality-hight-dict-image}
For any $q$,
$R^q\pi_*\underline{\mathbb{Q}}$ 
and
$R^q\pi_!\underline{\mathbb{Q}}$ 
are constant sheaves on $\chair$.
\end{lemma}

\proof
Since $\pi$ is locally trivial with respect to 
%the classical topology,
the complex topology,
$R^q\pi_*\underline{\mathbb{Q}}$ 
and
$R^q\pi_!\underline{\mathbb{Q}}$
are locally constant.
%Note that $\chair$ is path-connected.
Take a base point $x$ in $\chair$.
The fundamental group
$\pi_1(\chair,x)$ acts on the 
stalks $(R^q\pi_*\underline{\mathbb{Q}})_x\cong
H^q(\pi^{-1}(x);\mathbb{Q})$
and
$(R^q\pi_!\underline{\mathbb{Q}})_x\cong
H^q_c(\pi^{-1}(x);\mathbb{Q})$
through the map
$\pi_1(\chair,x)\to \pi_0(\pgl)$.
Since $\pgl$ is path-connected,
the action is trivial.
Hence the locally constant sheaves
$R^q\pi_*\underline{\mathbb{Q}}$ 
and
$R^q\pi_*\underline{\mathbb{Q}}$
are constant.  
\qed

\bigskip

By the above argument,
we can apply
%Theorem~\ref{thm:product-decomposition-VHP-general},
Theorem~\ref{thm:VHP-product-complex-topology}.
Then we obtain the following theorem.

%\begin{conjecture}\label{conj:VHPchtoair}\rm
\begin{theorem}[Theorem~\ref{thm:Main-Theorem2}]\label{thm:VHPchtoair}
The virtual Hodge polynomial of 
the character variety $\chtoair$
of absolutely irreducible representations
is given by 
\[ \begin{array}{rcl}
    {\rm VHP}(\chtoair)&=&\displaystyle
   \frac{(1-z^m)(1-z^{m-1})}{1-z^2},\\[3mm]
    {\rm VHP}_c(\chtoair)&=&\displaystyle
   \frac{z^{2m}(z^m-1)(z^{m-1}-1)}{z^2-1}.\\
   \end{array}\]
%\end{conjecture}
\end{theorem}

\proof
The theorem follows from
Theorem~\ref{thm:VHPreptoair}.
%Theorem~\ref{thm:product-decomposition-VHP-general}, 
%and
%Lemma~\ref{lemma:triviality-hight-dict-image}.
%Proposition~\ref{prop:Leray-product}.
%and
%Remark~\ref{remark:VHP_c-product-decomposition}.
\qed

\if0
\proof
The theorem may be proved as follows.
If $\pglto\to\reptoair\to\chtoair$ 
is a fiber bundle with respect to the Zariski topology,
then ${\rm VHP}_c(\reptoair)={\rm VHP}_c(\chtoair)\cdot
{\rm VHP}_c(\pglto)$.
Since ${\rm VHP}_c(\pglto)=z(z^2-1)$,
we may obtain that ${\rm VHP}_c(\chtoair)=
{\rm VHP}_c(\reptoair)/z(z^2-1)$.
\qed
\fi

%\input{counting}

\section{The number of absolutely irreducible representations}
\label{section:counting}

Let $p$ be a prime number and let $q$ be a power of $p$.
We denote by $\fq$ the finite field with $q$ elements.
In this section we study the number of
absolutely irreducible representations of degree $2$
over $\fq$
for the free monoid with $m$ generators.
For a scheme $X$, 
we denote by $|X(\fq)|$ the number of $\fq$-valued points of $X$.
Then $|{\rm Rep}_n(m)_{air}(\fq)|$
is the number of absolutely irreducible representations
of degree $n$ over $\fq$
for the free monoid with $m$ generators.
We show that 
$|{\rm Ch}_n(m)_{air}(\fq)|$ is the number of isomorphism classes of
such representations.
In case $n=2$,
we show that these numbers coincide with the virtual
Hodge polynomials evaluated at $q$.

Let $\fqbar$ be the algebraic closure of $\fq$ and
let $F: \fqbar\to \fqbar$ be the Frobenius map
given by $F(x)=x^q$. 
If a scheme $X$ is defined over $\mathbb{Z}$,
$F$ induces a map $F: X(\fqbar)\to X(\fqbar)$.
Since $F$ is a topological generator of ${\rm Gal}(\fqbar/\fq)$,
$|X(\fq)|$ is the number of fixed points under $F$.

\begin{proposition}\label{prop:number-chsc}
The number of $\fq$-valued points of 
${\rm Ch}_2(m)_{ss}$ is given by
\[ |{\rm Ch}_2(m)_{ss}(\fq)|=q^m(q^m-1).\]
\end{proposition}

\proof
By Corollary~\ref{cor:chss-unconfig},
we have an isomorphism
\[ {\rm Ch}_2(m)_{ss}(\fqbar)\cong \mathcal{C}_2(\fqbar{}^m). \]
Take $x=[a,b]\in \mathcal{C}_2(\fqbar{}^m)$.
Then $F(x)=[F(a),F(b)]$.
If $F(x)=x$,
then (i) $F(a)=a$ and $F(b)=b$,
or (ii) $F(a)=b$ and $F(b)=a$.
We denote by $X_1$ and $X_2$ the subsets of 
${\rm Ch}_2(m)_{ss}(\fq)$
consisting of elements
of type (i) and (ii), respectively.

In case (i),
$a, b\in (\fq)^m$.
Since $|F_2(\fq{}^m)|=q^m(q^m-1)$
and $\Sigma_2$ freely acts on $F_2(\fq{}^m)$, 
the number of $X_1$ is given by 
\[ |X_1|=\frac{1}{2}q^m(q^m-1). \]

In case (ii),
$a\in (\mathbb{F}_{q^2})^m$ and $b=F(a)$.
Hence $b$ is determined by $a$.
Since $a\neq b$,
$a\in (\mathbb{F}_{q^2})^m-(\fq)^m$.
Noticing $[a,F(a)]=[F(a),a]$,
we obtain that the number of $X_2$ is given by
\[ |X_2|=\frac{1}{2}(q^{2m}-q^m).\]

The number of $\fq$-valued points of ${\rm Ch}_2(m)_{ss}$ is
calculated as
\[ |X_1|+|X_2|=q^m(q^m-1).\]
\qed

%\begin{remark}\rm
%We can prove Proposition~\ref{prop:number-chsc}
%by using \'etale cohomology as follows.
%\end{remark}

\begin{proposition}\label{prop:number-repsc}
The number of $\fq$-valued points of
${\rm Rep}_2(m)_{ss}$ is given by
\[ |{\rm Rep}_2(m)_{ss}(\fq)|= q^{m+2}(q^m-1).\]
\end{proposition}

\proof
By Theorem~\ref{theorem:repss-description},
we have an isomorphism
\[ {\rm Rep}_2(m)_{ss}(\fqbar)\cong
{\rm PGL}_2(\fqbar)\times_{\hto} F_2(\fqbar{}^m). \]
Take $x=[G;a,b]\in {\rm Rep}_2(m)_{ss}(\fqbar)$,
where $G\in {\rm PGL}_2(\fqbar)$ and 
$(a,b)\in F_2(\fqbar{}^m)$. 
Let us regard $\tto_2$
as a closed subgroup scheme of ${\rm PGL}_2$.
If $F(x)=x$, then
(i) there exists $T\in {\rm T}_2(\fqbar)$
such that $F(G)=GT$ and $F(a)=a, F(b)=b$, or
(ii) there exists $T\in {\rm T}_2(\fqbar)$
such that $F(G)=GT\tau$
and $F(a)=b, F(b)=a$,
where $\tau$ is the permutation matrix
corresponding to the permutation $(1,2)$.
We denote by $X_1$ and $X_2$
the subsets of ${\rm Rep}_2(m)_{ss}(\fq)$
consisting of elements of type (i)
and (ii), respectively.

In case (i),
$(a,b)\in F_2(\fq{}^m)$.
We can take $T'\in {\rm T}_2(\fqbar)$
such that $T'F(T')^{-1}=T$.
Setting $G'=GT'$,
we have $F(G')=G'$, 
and hence
$G'\in {\rm PGL}_2(\fq)$.
Then
$x=[G; a,b]=[G'; a,b]$.
This means we can take a representative of $x$
as $[G';a,b]$ with
$G'\in {\rm PGL}_2(\fq)$ and
$(a,b)\in F_2(\fq{}^m)$.
Let $[\widetilde{G}; \widetilde{a},\widetilde{b}]$ 
be another representative of $x$
with $\widetilde{G}\in {\rm PGL}_2(\fq)$
and $(\widetilde{a},\widetilde{b})\in F_2(\fq{}^m)$.
Then there exists $H\in \hto$ such that
$\widetilde{G}=G'H$ and
$H^{-1}(\widetilde{a},\widetilde{b})H=(a,b)$.
In particular, 
$H={G'}^{-1}\widetilde{G}\in \hto\cap {\rm PGL}_2(\fq)$.
Since $\hto\cap {\rm PGL}_2(\fq)$  freely acts
on ${\rm PGL}_2(\fq)\times F_2(\fq{}^m)$,
the number of $X_1$ is given by

\[ \begin{array}{rcl}
     |X_1|&=&{\displaystyle \frac{|{\rm PGL}_2(\fq)|\cdot |F_2(\fq{}^m)|}
             {|\hto\cap {\rm PGL}_2(\fq)|}}\\[4mm]
          &=&{\displaystyle \frac{1}{2}q^{m+1}(q^m-1)(q+1)}.\\
   \end{array}\]

In case (ii),
$b=F(a)$ 
with $a\in (\mathbb{F}_{q^2})^m-(\fq)^m$.
We can take $T'\in {\rm T}_2(\fqbar)$
such that $T'F^2(T')^{-1}=T$.
Setting $G'=GT'\tau F(T')\tau$, 
we have $F(G')=G'\tau$.
In particular,
$F^2(G')=G'$ and $G'\in {\rm PGL}_2(\mathbb{F}_{q^2})$.
Notice that 
$G'$ has the form
$G'=\left(\begin{array}{cc}
                         1 & 1   \\
                         s & s^q \\
                        \end{array} \right)
    \left(\begin{array}{cc}
           t & 0        \\
           0      & t^q \\
          \end{array}\right)$
with $s\in \mathbb{F}_{q^2}-\fq,
t \in \mathbb{F}_{q^2}^{\times}$.
Put $S=\left(\begin{array}{cc}
                         1 & 1   \\
                         s & s^q \\
                        \end{array} \right)$.
Then $x=[S;a,F(a)]$.
So any element of $X_2$ has a representative 
of this form.
Since any element of $X_2$
has exactly two representatives of this form,
the number of $X_2$ is given by
\[ |X_2|=\frac{1}{2}q^{m+1}(q^m-1)(q-1).\]

The number of $\fq$-valued points of ${\rm Rep}_2(m)_{ss}$
is calculated as
\[ %\frac{1}{2}q^{m+1}(q^m-1)(q+1)+
   %\frac{1}{2}q^{m+1}(q^m-1)(q-1)=
   |X_1|+|X_2|=           q^{m+2}(q^m-1).\]
\qed

\begin{proposition}
\label{prop:fq-count-repnmsc}
The number of $\fq$-valued points of
${\rm Rep}_n(m)_{sc}$ and ${\rm Ch}_n(m)_{sc}$ 
are given by
\[  |{\rm Rep}_2(m)_{sc}(\fq)|= 
       |{\rm Ch}_2(m)_{sc}(\fq)|= q^m.\]
\end{proposition}

\proof
The proposition follows from Theorems~\ref{thm:descriptionRepsc}
and \ref{thm:descriptionChsc}.
\qed

\bigskip

In order to calculate the number of $\fq$-valued points
of the other moduli spaces, 
we need the following lemmas.

\begin{lemma}\label{lemma:key-lemma-counting}
Let $L$ be a field and let $\overline{L}$ be
its separable closure.
Let $G$ be an algebraic group over $L$
and let $X$ be a scheme of finite type over $L$.
Suppose that $G$ acts on $X$ over $L$.
If $G(\overline{L})$ freely acts on $X(\overline{L})$
and the Galois cohomology $H^1(L;G)$ is trivial,
then there is a bijection
\[ (X(\overline{L})/G(\overline{L}))^{{\rm Gal}(\overline{L}/L)}
   \cong X(L)/G(L).\]
\end{lemma}

\proof
We put ${\rm Gal}={\rm Gal}(\overline{L}/L)$.
Let $x$ be an element in $X(\overline{L})$ such that
the image under the quotient map
$X(\overline{L})\to X(\overline{L})/G(\overline{L})$
is invariant under the action of ${\rm Gal}$.
Since the action of $G(\overline{L})$ 
on $X(\overline{L})$ is free,
there exists a unique $c(\sigma)\in G(\overline{L})$
such that $x^{\sigma}=c(\sigma)x$ for each 
$\sigma\in {\rm Gal}$.
Then $c$ is a continuous
$1$-cocycle for ${\rm Gal}$
with the values in $G(\overline{L})$.
Since $H^1(L;G)$ is trivial,
there exists $g\in G(\overline{L})$ such that
$c(\sigma)=g^{\sigma}g^{-1}$ 
for all $\sigma\in {\rm Gal}$.
Then $(g^{-1}x)^{\sigma}=g^{-1}x$ and hence
$g^{-1}x\in X(L)$.
This means the canonical map
$X(L)\to (X(\overline{L})/G(\overline{L}))^{{\rm Gal}}$
is surjective.

Let $x_1$ and $x_2$ be elements in $X(L)$
that coincide in 
$(X(\overline{L})/G(\overline{L}))^{{\rm Gal}}$.
There exists $h\in G(\overline{L})$
such that $x_1=hx_2$.
Then $x_1=h^{\sigma}x_2$ for any $\sigma\in
{\rm Gal}$.
Since the action of $G(\overline{L})$ on $X(\overline{L})$
is free, $h^{\sigma}=h$ and hence $h\in G(L)$.
Therefore,
the map $X(L)\to 
(X(\overline{L})/G(\overline{L}))^{{\rm Gal}}$
induces a bijection 
\[ X(L)/G(L)\stackrel{\cong}{\longrightarrow}
   (X(\overline{L})/G(\overline{L}))^{{\rm Gal}}.\]
\qed

\begin{lemma}[{\cite{Lang}},
see also {\cite[Chapter~VI, Proposition~3]{Serre3}}]
%{see also \cite[Chapter~III, \S2.3, Theorem~1']{Serre2}}]
\label{lemma:general-vanishing-Lang}
The Galois cohomology $H^1(\fq;G)$ is trivial
for any connected algebraic group $G$ over $\fq$.
\end{lemma}

\begin{proposition}\label{prop:fq-number-rep2mu}
The number of $\fq$-valued points of ${\rm Rep}_2(m)_u$
is given by
\[ |{\rm Rep}_2(m)_u(\fq)|=q^m(q^m-1)(q+1).\]
\end{proposition}

\proof
By Theorem~\ref{theorem:description-repun},
the map
${\rm PGL}_2\times_{\bto}\utom\to {\rm Rep}_2(m)_u$
of algebraic varieties induces a bijection
\[ {\rm PGL}_2(\fqbar)\times_{\bto(\fqbar)}\utom(\fqbar)
    \stackrel{\cong}{\longrightarrow} 
   {\rm Rep}_2(m)_u(\fqbar). \]
Since this bijection is compatible
with the action of the Galois group ${\rm Gal}(\fqbar/\fq)$
on both sides,
we see that
\[ {\rm Rep}_2(m)_{u}(\fq)\cong
   ({\rm PGL}_2(\fqbar)\times_{\bto(\fqbar)}\utom(\fqbar))^{{\rm Gal}(\fqbar/\fq)}.\]
Notice that the action of $\bto(\fqbar)$ on 
${\rm PGL}_2(\fqbar)\times \utom(\fqbar)$
is free.
By
Lemmas~\ref{lemma:key-lemma-counting} 
%and \ref{lemma:H^1-trivial-B},
and \ref{lemma:general-vanishing-Lang},
we obtain
\[ |{\rm Rep}_2(m)_u(\fq)|=
    \frac{|{\rm PGL}_2(\fq)|\cdot|\utom(\fq)|}
         {|\bto(\fq)|}.\]
The proposition follows from the fact that
$\utom(\fq)\cong \fq{}^m\times (\fq{}^m-0)$.
\qed

\begin{proposition}\label{prop:number-chtomu}
The number of $\fq$-valued points of 
${\rm Ch}_2(m)_u$ is given by
\[ |{\rm Ch}_2(m)_u(\fq)|=\frac{q^m(q^m-1)}{q-1}.\]
\end{proposition}

\proof
\if0
By definition,
there is a bijection
\[ {\rm Ch}_2(m)_u(\fqbar)\cong 
    {\rm Rep}_2(m)_u(\fqbar)/{\rm PGL}_2(\fqbar)\]
This implies that there is a bijection
\[ {\rm Ch}_2(m)_u(\fqbar)\cong 
   \utom(\fqbar)/\tto_2(\fqbar), \]
where we regard $\tto_2$ as a closed
subgroup scheme of ${\rm PGL}_2$.
Notice that the action of $\tto_2(\fqbar)$
on $\utom(\fqbar)$ is free.
%Since $\tto_2\cong \mathbb{G}_m$,
Note that the Galois cohomology $H^1(\fq;\tto_2)$ is trivial
by Lemma~\ref{lemma:general-vanishing-Lang}.
By Lemma~\ref{lemma:key-lemma-counting},
we have
\[ |{\rm Ch}_2(m)_u(\fq)|=
   \frac{|\utom(\fq)|}{|\tto_2(\fq)|}.\]
The proposition is easily obtained from this. 
\fi
By Corollary~\ref{cor:descriptopn-chun},
the map
$\utom/\bto\to {\rm Ch}_2(m)_u$
of algebraic varieties induces a bijection 
\[ (\fqbar)^m\times {\rm P}^{m-1}(\fqbar)\stackrel{\cong}{\longrightarrow}
   {\rm Ch}_2(m)_u(\fqbar)\]
of $\fqbar$-valued points.
This bijection is compatible with the action of 
the Galois group ${\rm Gal}(\fqbar/\fq)$
on both sides, and
hence we obtain the following bijection
\[ (\fq)^m\times {\rm P}^{m-1}(\fq)\stackrel{\cong}{\longrightarrow}
   {\rm Ch}_2(m)_u(\fq).\]
The proposition follows from the fact that
$|(\fq)^m\times {\rm P}^{m-1}(\fq)|=q^m(q^m-1)/(q-1)$.
\qed

\if0
\proof
By definition,
there is a bijection
\[ {\rm Ch}_2(m)_u(\fqbar)\cong 
    {\rm Rep}_2(m)_u(\fqbar)/{\rm PGL}_2(\fqbar)\]
This implies that there is a bijection
\[ {\rm Ch}_2(m)_u(\fqbar)\cong 
   {\rm U}_2(m)(\fqbar)/{\rm T}_2(\fqbar). \]
Let $x\in {\rm Ch}_2(m)(\fq)$
and let $y\in \utom(\fqbar)$ be a 
representative of $x$.
Since $F(x)=x$,
there exists $T\in\tto_2(\fqbar)$
such that $F(y)=Ty$.
We can take $T'\in \tto_2(\fqbar)$
such that $F(T')T'^{-1}=T$.
Setting $y'={T'}^{-1}y$,
we have $F(y')=y'$
and hence $y'\in \utom(\fq)$.
Note that $y'$ is also a representative of $x$.

Let $y''\in \utom(\fq)$ be another
representative of $x$.
There exists $T''\in\tto_2(\fqbar)$
such that $y'=T''y''$.
Then $y'=F(T'')y''$ since
$F(y')=y'$ and $F(y'')=y''$.
Noticing that ${\rm T}_2(\fqbar)$ freely acts on
${\rm U}_2(m)(\fqbar)$,
we see that $F(T'')=T''$ and $T''\in \tto_2(\fq)$.
We obtain that
\[ {\rm Ch}_2(m)_u(\fq)\cong {\rm U}_2(m)(\fq)/{\rm T}_2(\fq),\]
and 
\[ \begin{array}{rcl}
    |{\rm Ch}_2(m)_u(\fq)|&=&
    {\displaystyle \frac{|{\rm U}_2(m)(\fq)|}
                   {|{\rm T}_2(\fq)|}}\\[4mm]
    &=&{\displaystyle \frac{q^m(q^m-1)}{q-1}}.\\
   \end{array}\]    
\qed
\fi

\begin{proposition}
\label{prop:counting-fq-B}
The number of $\fq$-valued points of ${\rm Rep}_n(m)_B$
is given by
\[ |{\rm Rep}_n(m)_B(\fq)|=\frac{
   q^{m(n-1)(n-2)/2}(q^m-q)^{n-1}
   \prod_{k=0}^{n-1}(q^m-k)\prod_{k=1}^n (q^k-1)}
   {(q-1)^n}.\]
\end{proposition}

\proof
By \cite[\S3]{topos},
we have a bijection
\[ {\rm Rep}_n(m)_B(\fqbar)\cong
   {\rm PGL}_n(\fqbar)\times_{{\rm B}_n(\fqbar)}
   {\rm B}_n(m)_B(\fqbar),\]
where we regard ${\rm B}_n(m)_B$ as a scheme over $\mathbb{Z}$.
Notice that the action of ${\rm B}_n(\fqbar)$
on ${\rm PGL}_n(\fqbar)\times {\rm B}_n(m)_B(\fqbar)$
is free.
By Lemmas~\ref{lemma:key-lemma-counting}
%and Lemma~\ref{lemma:H^1-trivial-B},
and \ref{lemma:general-vanishing-Lang},
we obtain
\[ |{\rm Rep}_n(m)_B(\fq)|=
   \frac{|{\rm PGL}_n(\fq)|\cdot |{\rm B}_n(m)_B(\fq)|}
        {|{\rm B}_n(\fq)|}\]
The proposition follows from the fact that
\[ |B_n(m)_B(\fq)|=q^{m(n-1)(n-2)/2}(q^m-q)^{n-1}
     \prod_{k=0}^{n-1}(q^m-k).\]
\qed

\if0
\proof
By \cite[\S3]{topos},
we have a bijection
\[ {\rm Rep}_n(m)_B(\fqbar)\cong
   {\rm PGL}_n(\fqbar)\times_{{\rm B}_n(\fqbar)}{\rm B}_n(m)_B(\fqbar).\]
In the same way as the proof of 
Proposition~\ref{prop:fq-number-rep2mu},
we obtain that
the number of $\fq$-valued points of 
${\rm Rep}_n(m)_B$ is
\[ \frac{|{\rm PGL}_n(\fq)| \cdot 
    |{\rm B}_n(m)_B(\fq)|}
    {|{\rm B}_n(\fq)|},\]
The proposition follows from the facts that
\[ |{\rm PGL}_n(\fq)/{\rm B}_n(\fq)|=
   \frac{(q^n-1)(q^{n-1}-1)\cdots (q-1)}
        {(q-1)^n},\]
and
\[ |B_n(m)_B(\fq)|=q^{m(n-1)(n-2)/2}(q^m-q)^{n-1}
     \prod_{k=0}^{n-1}(q^m-k).\]
\qed
\fi

\begin{proposition}
The number of $\fq$-valued points of
${\rm Ch}_n(m)_B$ is given by
\[ |{\rm Ch}_n(m)_B(\fq)|=
     \frac{(q^{m-1}-1)^{n-1}q^{(m-1)(n-1)(n-2)/2}
     \prod_{k=0}^{n-1}(q^m-k)}{(q-1)^{n-1}}.\] 
\end{proposition}

\proof
By definition,
there is a bijection
\[ {\rm Ch}_n(m)_B(\fqbar)\cong
     {\rm Rep}_n(m)_B(\fqbar)/{\rm PGL}_n(\fqbar)\]
This implies that 
there is a bijection
\[ {\rm Ch}_n(m)_B(\fqbar)\cong
     {{\rm B}_n(m)_B(\fqbar)/{\rm B}_n(\fqbar)}\]
Note that ${\rm B}_n(\fqbar)$ freely acts on
${\rm B}_2(m)_B(\fqbar)$.
By Lemmas~\ref{lemma:key-lemma-counting} 
%and \ref{lemma:H^1-trivial-B},
and \ref{lemma:general-vanishing-Lang},
we obtain
\[ |{\rm Ch}_n(m)_B(\fq)|=
   \frac{|{\rm B}_n(m)_B(\fq)|}
        {|{\rm B}_n(\fq)|}. \]
The proposition is easily obtained from this.
\qed

\if
\proof
By definition,
there is a bijection
\[ {\rm Ch}_n(m)_B(\fqbar)\cong
     {\rm Rep}_n(m)_B(\fqbar)/{\rm PGL}_n(\fqbar)\]
This implies that 
there is a bijection
\[ {\rm Ch}_n(m)_B(\fqbar)\cong
     {{\rm B}_n(m)_B(\fqbar)/{\rm B}_n(\fqbar)}\]
Note that ${\rm B}_n(\fqbar)$ freely acts on
${\rm B}_2(m)_B(\fqbar)$.
By Lemma~\ref{lemma:H^1-trivial-B},
$H^1({\rm Gal}(\fqbar/\fq);{\rm B}_n)$
is trivial.
In the same way as the proof of 
Proposition~\ref{prop:number-chtomu},
we obtain that
\[ {\rm Ch}_n(m)_B(\fq)\cong {\rm B}_n(m)_B(\fq)/{\rm B}_n(\fq)\]
and
\[ \begin{array}{rcl}
    |{\rm Ch}_n(m)_m(\fq)|&=&
    {\displaystyle \frac{|{\rm B}_n(m)_B(\fq)|}
                   {|{\rm B}_2(\fq)|}}\\[4mm]
    &=&{\displaystyle \frac{(q^{m-1}-1)^{n-1}q^{(m-1)(n-1)(n-2)/2}
     \prod_{k=0}^{n-1}(q^m-k)}{(q-1)^{n-1}}}.\\
   \end{array}\]    
\qed
\fi

\if0
(We have to check the conditions
to hold the following statement.)
If $X$ is defined over $\mathbb{Z}$,
then we have  
\[ |X(\fq)|={\rm VHP}_c(X_{\mathbb{C}})(\qh,\qh),\] 
where $X_{\mathbb{C}}$ is an algebraic scheme over $\mathbb{C}$
obtained by the base change.
(We have to check references)
By Theorem~\ref{thm:VHPreptoair} and \ref{thm:VHPchtoair},
we obtain the following theorems. 
\fi

\begin{theorem}[{cf.~Theorem~\ref{main-theorem-fq-counting}}]
\label{thm:number-fq-air}
The number of absolutely irreducible representations 
of degree $2$ over $\fq$
for the free monoid with $m$ generators 
coincides with
the virtual Hodge polynomial evaluated at $q$
so that
\[ \begin{array}{rcl}
     |{\rm Rep}_2(m)_{air}(\fq)|&=&
    {\rm VHP}_c({\rm Rep}_2(m)_{air}(\mathbb{C}))(q)\\[2mm]
    &=&q^{2m+1}(q^m-1)(q^{m-1}-1).\\
   \end{array}\] 
\end{theorem}

\proof
\if0
The scheme ${\rm Rep}_2(m)_{air}$ is an open subscheme
of ${\rm Rep}_2(m)={\rm M}_2(\mathbb{A})^m$.
Let $(A_1,\ldots,A_m)\in {\rm M}_2(\fqbar)^m$.
We denote by
$\fqbar\langle A_1,\ldots,A_m\rangle$ 
the $\fqbar$-subalgebra of ${\rm M}_n(\fqbar)$
generated by $A_1,\ldots,A_m$.
Then
$(A_1,\ldots,A_m)\in {\rm Rep}_2(m)_{air}(\fqbar)$
if and only if 
$\dim_{\fqbar}\fqbar\langle A_1,\ldots, A_m\rangle=4$.
We suppose that
$\dim_{\fqbar}\fqbar\langle A_1,\ldots, A_m\rangle<4$.
Then there exists $P\in {\rm PGL}_2(\fqbar)$
such that $P^{-1}A_iP$ is an upper triangular matrix
for $i=1,\ldots,m$, and
$(P;P^{-1}A_1P,\ldots,P^{-1}A_mP)\in {\rm Rep}_2(m)_{\ast}(\fqbar)$
for $\ast=sc,ss,u$ or $B$.
%Since $(A_1,\ldots,A_m)\in {\rm M}_2(\fq)^m$,
%we see that 
%$(P;P^{-1}A_1P,\ldots,P^{-1}A_mP)\in {\rm Rep}_2(m)_{\ast}(\fq)$.
Since the map
\[ \begin{array}{ccc}
   {\rm Rep}_2(m)_*(\fqbar)&\longrightarrow&
   {\rm Rep}_2(m)(\fqbar)={\rm M}_2(\fqbar)^m\\[2mm]
   (P;A_1,\ldots,A_m)&\mapsto& (PA_1P^{-1},\ldots,PA_mP^{-1})
  \end{array} \]
is injective for $*=sc,ss,u,B$ by construction,
we can regard ${\rm M}_2(\fqbar)^m$ 
as the disjoint union of 
${\rm Rep}_2(m)_*(\fqbar)$ for $*=sc,ss,u,B,air$:
\[ {\rm M}_2(\fqbar)^m=\coprod_{*=sc,ss,u,B,air} {\rm Rep}_2(m)_*(\fqbar).\]
\fi
By definition,
we have a decomposition
\[ {\rm M}_2(\fqbar)^m=\coprod_{*=sc,ss,u,B,air} {\rm
  Rep}_2(m)_*(\fqbar).\]
This decomposition is compatible with the map $F$.
Taking the fixed points of $F$,
we obtain a decomposition of ${\rm M}_2(\fq)^m$:
\[     {\rm M}_2(\fq)^m=\coprod_{*=sc,ss,u,B,air}{\rm
  Rep}_2(m)_*(\fq).\] 
Hence 
\[ |{\rm Rep}_2(m)_{air}(\fq)|=
    {\displaystyle q^{4m}-\sum_{*=sc,ss,u,B}|
     {\rm Rep}_2(m)_*(\fq)|}.\]
The theorem follows from
Propositions~\ref{prop:number-repsc}, 
\ref{prop:fq-count-repnmsc},
\ref{prop:fq-number-rep2mu}, 
and \ref{prop:counting-fq-B}.
\qed

\bigskip

\if0
In order to calculate 
$|{\rm Ch}_2(m)_{air}(\fq)|$
we need the following lemma.

\begin{lemma}\label{lemma:Hilbert_Thm90}
The Galois cohomology
$H^1(\fq;{\rm PGL}_n)$
is trivial.
\end{lemma}

\proof
We have
$H^1(\fq;{\rm GL}_n)=0$
(see, for example, \cite[Chapter~III, \S1.1, Lemma~1]{Serre2}
or \cite[Chapter~X, \S1, Proposition~3]{Serre1}),
and
$H^2(\fq;\mathbb{G}_m))=0$
(see, for example, \cite[Chapter~II, \S3]{Serre2}
or \cite[Chapter~X, \S7]{Serre1}).
The lemma follows from
the long exact sequence 
associated to the following exact sequence
\[ 1\longrightarrow \mathbb{G}_m
    \longrightarrow {\rm GL}_n
    \longrightarrow {\rm PGL}_n\longrightarrow 1\]
of group schemes.
\qed
\fi

\begin{theorem}[{cf.~Theorem~\ref{main-theorem-fq-counting}}]
\label{th:numberofchair}
The number of isomorphism classes of 
absolutely irreducible representations 
of degree $n$ over $\fq$
for the free monoid with $m$ generators is 
\[ |{\rm Ch}_n(m)_{air}(\fq)|=\frac{|{\rm Rep}_n(m)_{air}(\fq)|}
                                 {|{\rm PGL}_n(\fq)|}.\]
In case $n=2$, 
the number coincides with
the virtual Hodge polynomial evaluated at $q$
so that
\[ \begin{array}{rcl}
     |{\rm Ch}_2(m)_{air}(\fq)|&=&
    {\rm VHP}_c({\rm Ch}_2(m)_{air}(\mathbb{C}))(q)\\[2mm]
     &=&\displaystyle \frac{q^{2m}(q^m-1)(q^{m-1}-1)}{q^2-1}.\\
   \end{array}\] 
%\[ |{\rm Ch}_2(m)_{air}(\fq)|=
%   \frac{q^{2m}(q^m-1)(q^{m-1}-1)}{q^2-1}.\]
\end{theorem}

\proof
Since ${\rm Ch}_n(m)_{air}$
is the geometric quotient of ${\rm Rep}_n(m)_{air}$
by ${\rm PGL}_n$,
there is a bijection
\[ {\rm Ch}_n(m)_{air}(\fqbar)\cong 
    {\rm Rep}_n(m)_{air}(\fqbar)/{\rm PGL}_n(\fqbar).\]
Notice that ${\rm PGL}_n(\fqbar)$ freely acts on
${\rm Rep}_n(m)_{air}(\fqbar)$.
By Lemmas~\ref{lemma:key-lemma-counting} 
%and \ref{lemma:Hilbert_Thm90},
and \ref{lemma:general-vanishing-Lang},
we obtain 
\[ |{\rm Ch}_n(m)_{air}(\fq)|=
   \frac{|{\rm Rep}_n(m)_{air}(\fq)|}
        {|{\rm PGL}_n(\fq)|}. \]
This shows that the number of isomorphism classes
of absolutely irreducible representations over $\fq$
is $|{\rm Ch}_n(m)_{air}(\fq)|$.
When $n=2$,
we can calculate the number 
by
Theorem~\ref{thm:number-fq-air}.
\qed

\if0
\begin{remark}\rm 
Let $y \in \chair(\fq)$. 
Since $H^{1}(\fq; {\rm PGL}_n) = \{1 \}$ by Hilbert's Theorem 90,  
the fiber of $\pi : \repair(\fq) 
\to \chair(\fq)$ at $y$ has a $\fq$-rational point.  
Recall that $\pi$ is a ${\rm PGL}_n$-principal fibre bundle. 
Hence $| \pi^{-1}(y)(\fq) | = |{\rm PGL}_n(\fq)|$.  
\end{remark}
\fi

\begin{remark}\rm
Let $X$ be a separated scheme of finite type 
over $\mathbb{Z}$.
If there exists a polynomial $P_X(t)\in\mathbb{Z}[t]$
such that $|X(\fq)|=P_X(q)$
for all finite fields $\fq$,
then
${\rm VHP}_c(X)$ is a polynomial
of $z=xy$, 
and 
\[ {\rm VHP}_c(X)(z)=P_X(z). \]
See \cite[\S6]{HRV} for more details.
\end{remark}

%\begin{remark}\rm 
%By Theorem \ref{th:numberofchair},  
%the number of isomorphism classes of 
%$2$-dimensional absolutely irreducible modules over $\fq$ of 
%the free algebra $\fq \langle X_1, X_2, \ldots, X_m \rangle$ 
%is equal to $\frac{q^{2m}(q^m-1)(q^{m-1}-1)}{q^2-1}$. 
%\end{remark}

%\input{ind}
\begin{remark}\rm 
Let $S_m$ be the quiver with one vertex and $m$ edge loops.
The path algebra of $S_m$ over a field $k$
is the free algebra $k \langle X_1, X_2, \ldots, X_m \rangle$.
Hence the representations of the quiver $S_m$
are the same things as those of the free monoid with $m$ generators.  
Let ${\rm AIR}_{S_m}(n,q)$
be the number of isomorphism classes of
$n$-dimensional absolutely irreducible representations
of $S_m$ over $\fq$.  
By Theorem~\ref{th:numberofchair},  
we have
\[ {\rm AIR}_{S_m}(2,q)=\frac{q^{2m}(q^m-1)(q^{m-1}-1)}{q^2-1}.\] 
Let ${\rm AID}_{S_m}(n,q)$ be the number of 
isomorphism classes of $n$-dimensional absolutely 
indecomposable representations of $S_m$ over $\fq$.
Using \cite[Theorem~4.6]{Hua},
we can calculate ${\rm AID}_{S_m}(2,q)$ as  
\[ {\rm AID}_{S_m}(2,q)=\frac{q^{2m-1}(q^{2m}-1)}{q^2-1}.\]
We can verify that this number is equal to the sum
\[ |{\rm Ch}_2(m)_u(\fq)|+|{\rm Ch}_2(m)_B(\fq)|
  +|{\rm Ch}_2(m)_{air}(\fq)|. \]
\end{remark}

\bigskip 

%\begin{remark}\rm 
Let $A_n(m)$ be the affine ring of ${\rm Rep}_n(m)$.  
Let $A_n(m)^{{\rm PGL}_n}$ be
the ${\rm PGL}_n$-invariant ring of $A_n(m)$. 
We set ${\rm Ch}_n(m) := {\rm Spec} A_n(m)^{{\rm PGL}_n}$. 
By \cite[Theorem~3]{Seshadri},
the set of $\fqbar$-valued points of
${\rm Ch}_n(m)$ consists of
the closed ${\rm PGL}_n$-orbits in ${\rm Rep}_n(m)$.
In particular, when $n=2$,
we have a decomposition
\[ {\rm Ch}_2(m)(\fqbar)={\rm Ch}_2(m)_{air}(\fqbar)\coprod 
                         {\rm Ch}_2(m)_{ss}(\fqbar)\coprod
                         {\rm Ch}_2(m)_{sc}(\fqbar)\]
of $\fqbar$-valued points.
This implies the following main theorem: 

\begin{theorem}[Theorem~\ref{theorem:Main-Theorem4}]\label{thm:VHPch} 
The number of ${\Bbb F}_q$-valued points of ${\rm Ch}_2(m)$ is given by 
\[ 
\begin{array}{ccl} 
|{\rm Ch}_2(m)({\mathbb F}_q)| & = & |{\rm Ch}_2(m)_{\rm air}({\mathbb F}_q)| 
+ |{\rm Ch}_2(m)_{ss}({\mathbb F}_q)| 
+ |{\rm Ch}_2(m)_{sc}({\mathbb F}_q)| \vspace*{1ex} \\ 
 & = & \displaystyle \frac{q^{2m+2}(q^{2m-3} - q^{m-2} - q^{m-3} + 1)}{q^2 - 1}.   
\end{array} 
\] 
In particular, the virtual Hodge polynomial of 
${\rm Ch}_2(m)$ is given by 
\[ 
%\begin{array}{ccl} 
{\rm VHP}_{c}({\rm Ch}_2(m))(z) 
  =  \displaystyle \frac{z^{2m+2}(z^{2m-3} - z^{m-2} - z^{m-3} + 1)}{z^2 - 1}.   
%\end{array} 
\]
\end{theorem} 
%\end{remark} 

\begin{remark}\rm
The Weil zeta functions of $\reptoairab, \chtoairab$, 
and ${\rm Ch}_2(m)$ are given by 
\[
\begin{array}{ccccc} 
Z(\reptoairab, q, t) & := & \displaystyle 
\exp\left(\sum_{n=1}^{\infty} 
\frac{| \reptoairab( \mathbb{F}_{q^n} ) |}{n} t^n \right) 
& = & \displaystyle 
\frac{(1-q^{3m+1}t)(1-q^{3m}t)}{(1-q^{4m}t)(1-q^{2m+1}t)},  \\[7mm]
Z(\chtoairab, q, t) & := & \displaystyle 
\exp\left(\sum_{n=1}^{\infty} 
\frac{| \chtoairab( \mathbb{F}_{q^n} ) |}{n} t^n \right) 
& = & \displaystyle  
\frac{\displaystyle  \prod_{i=1}^{[ \frac{m}{2} ]} (1-q^{2m+2i-2}t)}%
{\displaystyle  \prod_{i=1}^{[ \frac{m}{2} ]} (1-q^{4m-2i-1}t)},  \\[15mm]
Z({\rm Ch}_2(m), q, t) & := & 
\displaystyle 
\exp\left(\sum_{n=1}^{\infty} 
\frac{| {\rm Ch}_2(m)( \mathbb{F}_{q^n} ) |}{n} t^n \right) 
& = & 
\displaystyle
\frac{Z(\chtoairab, q, t)}{1-q^{2m}t} .  \\[5mm]
\end{array}
\]
The Hasse-Weil zeta functions of $\reptoairab, \chtoairab$,
and ${\rm Ch}_2(m)$ are given by 
\[
\begin{array}{ccccc} 
\zeta(\reptoairab, s) & := & \displaystyle \prod_{p} 
Z(\reptoairab, p, p^{-s}) 
& = & \displaystyle 
\frac{\zeta(s-4m) \zeta(s-2m-1)}{\zeta(s-3m-1) \zeta(s-3m)},  \\[5mm]
\zeta(\chtoairab, s) & := & \displaystyle \prod_{p} 
Z(\chtoairab, p, p^{-s}) 
& = & \displaystyle  
\frac{\displaystyle  \prod_{i=1}^{[ \frac{m}{2} ]} \zeta(s-4m+2i+1)}%
{\displaystyle  \prod_{i=1}^{[ \frac{m}{2} ]} \zeta(s-2m-2i+2)},  \\[12mm]
\zeta({\rm Ch}_2(m), s) & := & 
\displaystyle \prod_{p} Z({\rm Ch}_2(m), p, p^{-s}) 
& = & \zeta(\chtoairab, s)
\zeta(s-2m),  \\[5mm]
\end{array}
\]
where $\zeta(s)$ is the Riemann zeta function. 
The completions of these zeta functions are defined as 
\[
\begin{array}{rcl} 
\hat{\zeta}(\reptoairab, s) & := & \displaystyle 
\frac{\hat{\zeta}(s-4m) 
\hat{\zeta}(s-2m-1)}{\hat{\zeta}(s-3m-1) \hat{\zeta}(s-3m)},  \\[5mm]
\hat{\zeta}(\chtoairab, s) & := & \displaystyle  
\frac{\displaystyle  \prod_{i=1}^{[ \frac{m}{2} ]} \hat{\zeta}(s-4m+2i+1)}%
{\displaystyle  \prod_{i=1}^{[ \frac{m}{2} ]} \hat{\zeta}(s-2m-2i+2)},  \\%[12mm]
%\hat{\zeta}({\rm Ch}_2(m), s) & := & 
%\hat{\zeta}(\chtoairab, s)
%\hat{\zeta}(s-2m),
\end{array}
\]
where $\hat{\zeta}(s) := \pi^{-s/2} \Gamma(s/2) \zeta(s)$ 
is the completion of the Riemann zeta function.  
Since $\hat{\zeta}(1-s) = \hat{\zeta}(s)$, 
we have the following functional equations 
\[
\begin{array}{rcl} 
\hat{\zeta}(\reptoairab, 6m+2-s) & = & 
\hat{\zeta}(\reptoairab, s),  \\[2mm] 
%\hat{\zeta}(\chtoairab, 6m-2-s) & = & 1/\hat{\zeta}(\chtoairab, s).  
\hat{\zeta}(\chtoairab, 6m-2-s) & = & \hat{\zeta}(\chtoairab, s)^{-1}.\\[2mm]
%\hat{\zeta}({\rm Ch}_2(m),??) &=& 
%\hat{\zeta}({\rm Ch_2}(m),s)??.\\  
\end{array} 
\]
\end{remark}

\bigskip

\if0
\begin{itemize}

\item The argument in the proof of Conjecture~\ref{conj:VHPchtoair}
is not true.

%\item In descriptions of moduli spaces,
%an isomorphism of complex manifolds  
%may be insufficient to use the calculation 
%of virtual Hodge polynomials.

\end{itemize}
\fi

\bigskip

\end{document}